\documentclass[12pt]{article}

\usepackage[utf8]{inputenc}
\usepackage{amsmath,IEEEtrantools}

\usepackage{color}
\usepackage{xcolor}
\usepackage{algorithm2e}
\usepackage{color}
\usepackage{graphicx}
\usepackage{amssymb}
\usepackage{amsmath}
\usepackage{amsfonts}
\usepackage{mathtools}
\usepackage{float}
\usepackage{algorithm2e}
\usepackage{graphicx}
\usepackage{caption}
\usepackage{subcaption}
\usepackage{multirow}
\usepackage{verbatim}
\usepackage{enumitem}

\newtheorem{rem}{Remark}

\let\realverbatim=\verbatim
\let\realendverbatim=\endverbatim
\renewcommand\verbatim{\par\addvspace{6pt plus 2pt minus 1pt}\realverbatim}
\renewcommand\endverbatim{\realendverbatim\addvspace{6pt plus 2pt minus 1pt}}
\newcommand{\ma}{\textcolor{black}}

\newcommand{\re}{\textcolor{red}}

\newcommand{\bse}{\begin{subequations}}
\newcommand{\ese}{\end{subequations}}
\newcommand{\vect}[1]{\mbox{\boldmath $#1$}}
\newcommand{\bes} {\begin{eqnarray*}}
\newcommand{\ees} {\end{eqnarray*}}
\newcommand{\ds}{\displaystyle}
\newcommand{\be}{\begin{eqnarray}}
\newcommand{\ee}{\end{eqnarray}}
\newcommand{\non} \nonumber

\newcommand{\ra}{\rightarrow}

\newcommand{\ex}{\vect{e}_x}
\newcommand{\er}{\vect{e}_r}

\newcommand{\upz}{u_{p,0}}
\newcommand{\unz}{u_{n,0}}
\newcommand{\uno}{u_{n,1}}
\newcommand{\upo}{u_{p,1}}

\newcommand{\wno}{w_{n,1}}
\newcommand{\wpo}{w_{p,1}}

\newcommand{\gil}{G_1^{(L)}}
\newcommand{\gis}{G_1^{(S)}}
\newcommand{\gsm}{G_1^{({\rm smooth})}}

\newcommand{\z}{\zeta}

\newcommand{\psdz}{\psi^{(d)}_0}
\newcommand{\psdo}{\psi^{(d)}_1}

\newcommand{\ie}{{\it i.e.}\ }
\newcommand{\eg}{{\it e.g.}\ }
\newcommand{\etal}{{\it et al.} \ }

\newcommand{\xx}{\vect{x}}
\newcommand{\J}{\vect{{\cal F}}}
\newcommand{\vep}{\varepsilon}
\newcommand{\ns}{n}
\newcommand{\ps}{p}
\newcommand{\pb}{\bar{P}}
\newcommand{\nb}{\bar{N}}
\newcommand{\jpb}{\bar{J}_p}
\newcommand{\jnb}{\bar{J}_n}

\newcommand{\phit}{\tilde{\Phi}}
\newcommand{\Psh}{\hat{\Psi}}

\newsavebox{\astrutbox}
\sbox{\astrutbox}{\rule[-5pt]{0pt}{20pt}}

\newcommand\dd{\ensuremath{\partial}}

\title{Asymptotic models for transport in large aspect ratio nanopores}

\author{B. Matejczyk\thanks{Department of Mathematics, University of Warwick, CV4 7AL Coventry, United Kingdom}, 
J.-F. Pietschmann\thanks{Institute for Computational and Applied Mathematics, WWU M\"unster, M\"unster 48149, Germany and Osnabr\"uck University, Institute of Mathematics, 49069 Osnabr\"uck, Germany} ,
M.-T. Wolfram\thanks{Radon Institute for Computational and Applied Mathematics, Austrian Academy of Sciences, Altenberger Strasse 69, 4040 Linz, Austria}$\;^{,*}$, G. Richardson\thanks{School of Mathematics, University of Southampton, Southampton, UK, SO17 1BJ}}
\begin{document}

\maketitle
\begin{abstract}
Ion flow in charged nanopores is strongly influenced by  the ratio of the Debye length to the pore radius.
We investigate the asymptotic behaviour of solutions to the Poisson-Nernst-Planck (PNP) system in narrow pore like
geometries and study the influence of the pore geometry and surface charge on ion transport. The physical properties of real pores motivate 
the investigation of distinguished asymptotic limits, in which either the Debye length and pore radius are comparable or the pore length 
is very much greater than its radius. 
This results in a Quasi-1D PNP model which can be further simplified, in the physically relevant limit of strong pore wall surface charge, to a fully one-dimensional model. Favourable comparison is made 
to the two-dimensional PNP equations in typical pore geometries. It is also shown that, for physically realistic parameters, the standard 1D Area Averaged PNP model for ion flow through a pore is a very poor approximation to the (real) two-dimensional solution to the PNP equations. This leads us to propose that the Quasi-1D PNP model derived here, whose computational cost is significantly less than two-dimensional solution of the PNP equations, should replace the use of the 1D Area Averaged PNP equations as a tool to investigate ion and current flows in ion pores. 
\end{abstract}

%%%%%%%%%%%%%%%%% Problem statement
\section{Introduction}

Solid-state nanopores are nanoscale holes in synthetic materials, such as silicon nitrite, graphene or polyethylene terepthalate (PET). They can be produced in a 
variety of lengths and shapes, with diameters ranging from a few nanometer to openings at the micrometer scale. There has been a tremendous increase in research on nanopores 
over the last decades, mostly initiated by their use as sensors for DNA or other biomolecules. In a typical experiment, one or more nanopores are placed into a bath containing an ionic solution with potentially different ionic concentrations on each side of the pore. Then, an additional external potential is applied and the current generated by the ions moving through the pore is measured. 
When acting as sensor, the current fluctuation as an unknown molecule traverses through the pore, allows the determination of its structural or chemical properties, \cite{Schneider2012}. Another interesting characteristic of many pores is their rectification behaviour, meaning that the measured current for positive applied external voltage differs from that for negative voltage. This effect has been extensively studied and is believed to originate from the  combination of geometric and electrostatic effects.
 Hence,  the pores are effectively acting as  diodes which make them useful as a building block for more complex circuits.\\
 In terms of mathematical modelling, the Poisson-Nernst-Planck (PNP) equations have been used successfully to describe the flow of ions 
through pores, \cite{vlassiouk2008nanofluidic,siwy_ss_2003,PhysRevE.76.041202}. They consist of a set of drift--diffusion equations for the density of ions, 
self consistently coupled to a Poisson equation to account for the electrostatic interactions due to the charge of the ions themselves and other charges present in the system. 
They were first introduced to describe carrier transport in semiconductor devices where they are known as drift-diffusion equations (DDE), 
\cite{markowich1986uniform,markowich1990semiconductor}, but have also been applied to many other systems, \eg batteries \cite{newman12} and
ion channels \cite{corry2000tests,horng2012pnp,Yang2005,Qiong2011,matejczyk2017multiscale}. 
The analysis of the DDE is well understood, see \cite{Mar,markowich1990semiconductor} and its asymptotic behaviour has been analysed in a variety of contexts, 
such as semiconductor devices \cite{markowich1986asymptotic, markowich1986uniform}, solar cells \cite{courtier17,foster13,foster14,richardson16,richardson17}, cell membrane action 
potentials \cite{george15,richardson2009multiscale} and batteries \cite{lai2011mathematical,richardson12}. Perhaps more importantly, these asymptotic methods provide techniques which allow to systematically derive simplified models from numerically challenging PDE models of the underlying physics. For example equivalent circuit models of solar cells have been derived from DDE models in \cite{foster14,foster13,richardson12opv}, surface polarization models for both action potentials in \cite{richardson2009multiscale} and for perovskite solar cells in \cite{courtier17,richardson16,richardson17}. Also effective medium models for lithium-ion batteries \cite{lai2011mathematical,richardson12} and organic photovoltaic devices \cite{richardson17} have been derived from DDE models. In a similar vein one of the major aims of this work, is to systematically apply asymptotic methods  to a DDE model of an electrolyte in a nanopore to derive a simplified, computationally tractable model of ion transport through the pore.

With increasing miniaturisation of semiconductor devices as well as the application to nanoscale systems, 
the influence of finite size effects on the transport behaviour became more important. However the classical PNP equations treats particles as point charges, without size, and omits 
particle-particle interactions such as volume exclusion effects. For this reason several extensions were introduced, for example by including finite volume effects already in a microscopic model
\cite{burger2012nonlinear} or by using density functional theory to account for quantum mechanical effects \cite{Gillespie2002}. Dreyer and co-workers proposed a thermodynamically consistent coupling to the
Navier-Stokes equations, which includes the velocity of the solvent in \cite{Dreyer2013}.\\

In this paper we study  the classical PNP model and analyse  the behaviour of its solutions in the case of the nanopores with  a very large aspect ratio, i.e 
the pore radius is much smaller than the length of the pore. Our analysis is motivated by the geometry and structure of typical polyethylene terepthalate (PET) nanopores,
which are  produced by irradiating a $12\mu$m thick PET foil with heavy ions and subsequent chemical etching, \cite{siwy_ss_2003}. The resulting pores are radially symmetric and the etching creates carboxyl groups at the pore walls at an estimated density of $1$ electron per nm$^2$. They have typical opening diameters of $4-200$nm and $200-1000$nm at their respective ends. Thus their aspect ratio is in the range of $0.0003-0.08$.
While the PNP equations are in principle still valid in this regime, from a computational point of view it is very expensive, if at all possible, to simulate such a pore completely. Furthermore, the rectification behaviour,  is experimentally determined from so-called I-V-curves which are obtained by measuring the current over a certain range of applied voltages. Hence, the PNP equations have to be solved several times, which makes the problem of computational cost even more important.\\
To overcome these issues we use tools from asymptotic analysis to derive a one-dimensional approximation of the full PNP system that is still able 
to capture the physical behaviour of the pore and can also be used for numerical computations. In a similar context a one-dimensional area averaged asymptotic model has previously been used to study  ion flow through  biological ion channels \cite{singer_gillespie_norbury_eisenberg_2008,singer2009poisson,ionreport}. However, as far as we are aware, there has been no direct comparison between numerical solutions to the full PNP equations, in appropriate geometries, and this 1d-model. Indeed there are good reasons to suppose, as has been pointed out by Chen \etal \cite{Chen2014}, that even the full 3D-PNP is incapable of adequately describing the behaviour of interactions between ions occurring on the atomic lengthscale in an ion channel. Chen \etal \cite{Chen2014} instead make use of an approach based on the Fokker-Planck equations which allow them to capture the important effects of direct inter-ion interactions in the narrow neck of the channel and can be shown to lead, via an asymptotic approximation, to a Markovian transition rate model of a type which is often used to provide a phenomenological description of ion-channel behaviour (see for example \cite{Ball2002}). The present work, however, is 
concerned with ion transport in nanopores which are considerably larger structures than ion channels and for which PNP type models provide an appropriate description.

We start by discussing the PNP equations as well as respective physical parameter ranges in Section \ref{s:pnp}.
In Section  \ref{medium_regime} we present an asymptotic solution to the model for pores with large aspect ratio and a radius comparable to 
the Debye length of the electrolyte. This asymptotic solution allows us to characterise the behaviour of the pore in terms of the solution to a one-dimensional model. 
In Section \ref{numerics} we compare results obtained from the one-dimensional asymptotic model to numerical solutions of the full equations in  axisymetric large aspect ratio  pores.

\section{The PNP equations}\label{s:pnp}

\begin{figure}[H]
     \begin{center}
       \includegraphics[width=0.9\textwidth]{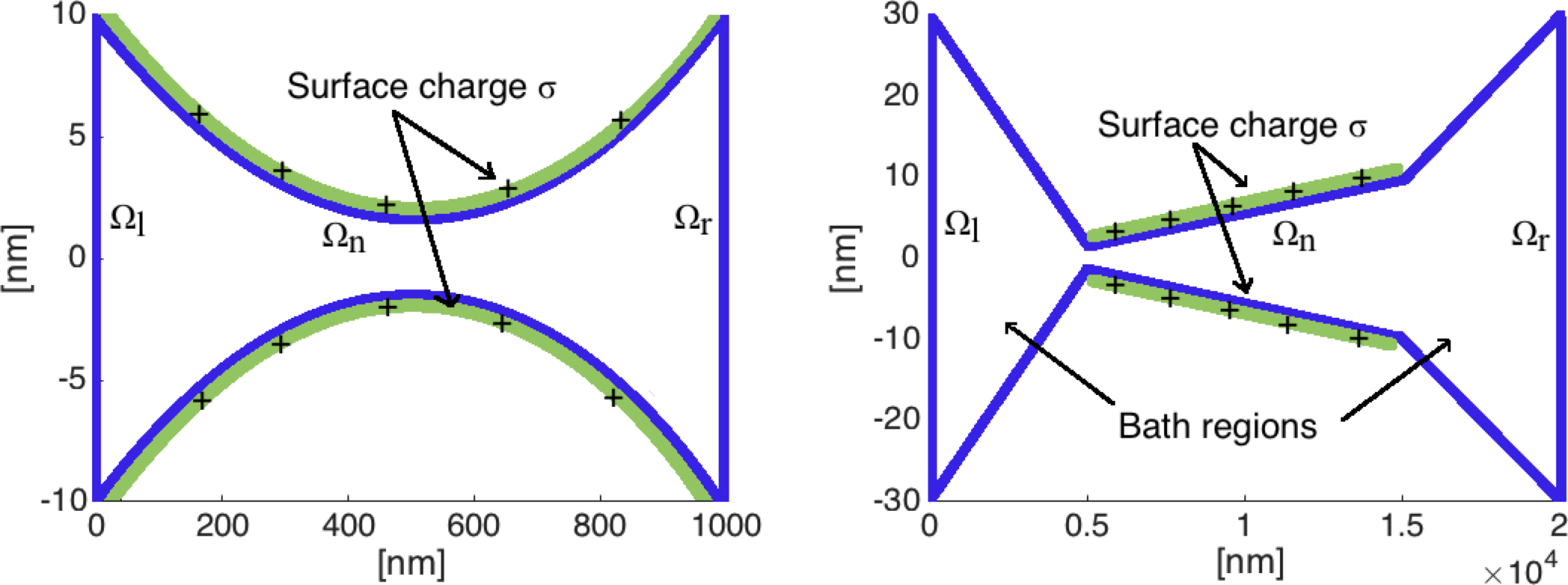}
        \caption{Sketches of the geometries considered for the nanopores.}  \label{f:geom}
       \label{geometries}
        \end{center}
\end{figure}

We start by presenting the mathematical model and its scaling which serves as the basis of our asymptotic analysis.
For ease of presentation, and because this is a typical set-up in practice, we restrict our attention to an ideal 1:1 electrolyte comprised of positive and negative ions of valency one and with concentrations $p^*$ and $n^*$ respectively (measured in moles per unit volume). 
Note that we use $*$ to indicate dimensional variables throughout the manuscript. \\
The PNP equations for the concentrations $n^* = n^*(\xx^*,t^*)$, $ p^* = p^*(\xx^*,t^*)$ and the electric potential $V^*= V^*(\xx^*,t^*) $ read as  \begin{IEEEeqnarray}{CC}
\IEEEyesnumber
\label{gen_PNP} \IEEEyessubnumber*
-\nabla^* \cdot ( \varepsilon  \nabla^* V^* ) = F(p^*-n^*),   \label{3dpoiss} \\
 \frac{\dd p^*}{\dd t^*}+ \nabla^* \cdot \J^*_p=0 , \qquad  \frac{\dd n^*}{\dd t^*}+ \nabla \cdot   \J_n^*=0 , \\ 
\J^*_n= -D_n \left(\nabla^* n^* - \frac{1}{V_T} n^* \nabla^* V^* \right)  \label{3dbv},\\
 \J^*_p = - D_p \left(\nabla^* p^* + \frac{1}{V_T}  p^* \nabla^* V^* \right).
\end{IEEEeqnarray} 
Here  $\J_p^*$ and $\J_n^*$  are the flux of positive and negative ions, respectively, $F$ is Faraday's constant, and $V_T$ the thermal voltage.   
The parameters $D_p$ and $D_n$ are the diffusion coefficients of the positive and negative ions, respectively, and the 
domain $\Omega$ is assumed axially symmetric being given by 
\begin{equation*}
\Omega= \{(x^*,y^*,z^*) : 0\leq x^* \leq L^*, \ma{0 \leq \sqrt{{y^*}^2+{z^*}^2} \leq  R^*(x^*)}\},
\end{equation*}
where $R^*(x^*)$ is the radius of the pore as a function of $x^*$.
The boundary of $\Omega$ is split into three subdomains, the left and the right entrance of the nanopore 
\begin{align*}
 \Omega_l = \{(x^*,y^*,z^*) \in \partial \Omega , x^*=0\} \text{ and }
 \Omega_r =\{(x^*,y^*,z^*) \in \partial \Omega , x^*=L\} ,
\end{align*}
as well as  the nanopore walls 
$
\Omega_N= \{(x^*,y^*,z^*) \in \partial \Omega , \sqrt{{y^*}^2+{z^*}^2}=\ma{R^*(x^*)}\}$. The considered geometry of the pore is depicted in Figure \ref{f:geom}. 
In the same Figure, we also present a more realistic  geometry, in which additional \emph{bath regions}  are attached at each end of the pore.
System \eqref{gen_PNP} is supplemented with the following boundary conditions:
\begin{align}
V^*|_{\Omega_l }= 0, \quad &V^*|_{ \Omega_r}= V_{appl}, \quad \{n^*,p^*\}|_{\Omega_l }= \{n_l,p_l\}, \quad \{n^*,p^*\}|_{ \Omega_r}= \{n_r,p_r\},  \label{Dir_BC}\\
&\J^*_p \cdot \vect{N}|_{\Omega_N} =\J^*_n \cdot \vect{N}|_{\Omega_N} =\vect{0},  \quad \left. \frac{\dd V^*}{\dd N^* } \right|_{\Omega_N }=\frac{\sigma^*(x^*) }{\vep}\label{Neu_BC}.
\end{align}
where $\dd/\dd N^*$ denotes the normal derivative to the pore boundary with respect to its unit outward normal $\vect{N}$, defined by
\bes
\vect{N}=\left( \er - \frac{d R^*}{d x^*} \ex \right) \left(1+ \left( \frac{d R^*}{d x^*}\right)^2\right)^{-1/2},
\ees
$\sigma(x^*)$ is the surface charge density on the pore wall and $\vep$ the permittivity of the electrolyte. The Dirichlet conditions \eqref{Dir_BC} correspond to a prescribed applied voltage and prescribed ion concentrations at each opening of the pore and in the bath regions, respectively. Here, for computational convenience, these are imposed on a fixed external boundary whereas it could be argued that these ought to be imposed as far-field conditions.
However these two sets of boundary conditions have almost identical solutions provided the pore is sufficiently wide when it is terminated by the artificial boundaries $\Omega_l$ and $\Omega_r$. 
Condition \eqref{Neu_BC} ensures that there is no ion flux through the pore walls and prescribes the fixed surface charge at these walls. We note that surface charge condition is asymptotically correct only if the permittivity $\vep$ of the 
electrolyte is much greater than that of the pore walls (which is the case for aqueous electrolytes); for more details see \cite{ionreport}.

The current-voltage curve (IV curve in short) is commonly used to characterise the behaviour of ion channels and nanopores. The respective current flow $I^*(x^*,t^*)$ can be computed by calculating the current flow through a cross-section on the pore, at $x^*=X^*$ say, being given by
\begin{align}
I^*(x^*,t^*)=F\hspace*{-4em} \int\limits_{ \{{y^*}^2+{z^*}^2 \leq {R^*}^2(X^*)\}\, \cap \,  \{x^*=X^*\} }\hspace*{-4em} \ex \cdot (\J^*_p -\J^*_n )|_{x^*=X^*} dS^*. \label{current}
\end{align}

\subsection{The 1D Area Averaged PNP equations}
The 1D Area Averaged PNP equations are a common reduction of the full PNP system, which is frequently used to calculate ion flux through radially symmetric nanopores because of its much reduced computational cost, see for example \cite{cervera_jcp_2006, cervera_epl_2005}. In this approach, at each point of the $x^*$ axis, the average of the ionic concentrations and the electrostatic potential on the disc $0 \leq \sqrt{{y^*}^2+{z^*}^2} \leq  R^*(x^*)$ is calculated. This yields \\
\begin{subequations}
\begin{align}
&\frac{\partial}{\partial x^*}\left(\varepsilon A(x^*) \frac{\partial V^*}{\partial x^*} \right) = -\left(F A(x^*) (p^*-n^*) +  \partial A(x^*)\sigma^*(x^*) \right), \label{1D_PNPa}\\
&D_n\frac{\partial}{\partial x^*} \left(A(x^*) \left(\frac{\partial n^*}{\partial x^*} - \frac{1}{V_T} n^* \frac{\partial V^*}{\partial x^*} \right) \right) = 0,\\
&D_p\frac{\partial}{\partial x^*} \left(A(x^*)\left(\frac{\partial p^*}{\partial x^*}  + \frac{1}{V_T} p^* \frac{\partial V^*}{\partial x^*} \right) \right)  = 0, \label{1D_PNPc}
\end{align}
\end{subequations}
where 
\begin{align}\label{area}
A(x^*) = R^*(x^*)^2 \pi
\end{align}
denotes the area and $\partial A(x^*) = 2 R(x^*) \pi$ the circumference of the pore at point $x^*$. 
The 1D Area Averaged PNP equations are based on the assumption that the influence of the surface charge on the ion concentration and the voltage can be averaged over the cross section of the pore. This assumption is only valid if the Debye length of the electrolyte is much greater than the pore width (e.g. for extremely dilute solutions) and does not hold for typical experimental conditions. 
%\re{An alternative 1D reduction proposed by Gillespie et al.{\cite{Gillespie2002}} for ion channels, in which the ion concentration and voltage are averaged over the equi-potential and equi-conentration surfaces. These surfaces correspond to the cross section of the ion channel in the case of radially symmetric distributed charges. However this assumption does not hold for nanopores and the acutual computation of these surfaces would be rather involved. \mage{(SEEMS TO BE SAME AS 1D AREA AVERAGED NONLINEAR PNP).} }}
Indeed, we shall show that solutions to the 1D  Area Averaged PNP equations do not provide a good approximation to the full 2D equations for the geometries we consider except in exceptional circumstances (e.g for extremely dilute electrolytes or very narrow pores) in which the Debye length of the electrolyte is very much greater than the pore width. The main aim of this work is to derive a similar 1D model that is capable of adequately approximating the 2D PNP equations. Comparison is made between the resulting 1D models and the full 2D system in Section \ref{numerics}.

%%%%%%Scaling 
\subsection{Scaling} \label{general_scaling_PNP}
We nondimensionalise system \eqref{gen_PNP} by introducing a  typical lateral lengthscale $L$,  a typical pore radius $R_0$, a typical concentration $ \bar c$ 
(measured in ion number per unit volume),  and a typical surface charge $\bar \sigma$. 
The great disparity in size between the lateral lengthscale $L$ and the pore radius $R_0$ motivates us to rescale differently in the these two dimensions.
This results in different scalings for the fluxes in the radial and lateral directions. 
We introduce the radial variables $r^*=\sqrt{{y^*}^2+{z^*}^2}$  and nondimensionalise as follows
\begin{align*}
& x^* = L x, ~~r^* = R_0 r, ~~p^* = \bar{c} p,~~ n^* = \bar{c} n,~~\sigma^* = \bar{\sigma}\sigma,~~V^* = V_T \phi\\
& \ds  \J_p^* \cdot \ex=\frac{\bar{D} \bar{c}}{L} u_p, ~~ \ds  \J_p^* \cdot \er= \frac{\bar{D} \bar{c} R_0}{L^2}w_p, ~~\ds  \J_n^* \cdot \ex=\frac{\bar{D} \bar{c}}{L} u_n, \ds  \J_n^* \cdot \er=\frac{\bar{D} \bar{c} R_0}{L^2}w_n,
\end{align*}
where $\bar{D}$ is a typical ionic diffusivity which we assume to be constant everywhere inside the domain. This leads to the following dimensionless formulation of system \eqref{gen_PNP}
\begin{subequations}
\begin{align}
&\frac{\dd p}{\dd t}+\frac{\dd u_p}{\dd x}+ \frac{1}{r} \frac{\dd}{\dd r} ( r w_p) =0, \label{nd1}\\
&\frac{\dd n}{\dd t}+\frac{\dd u_n}{\dd x}+ \frac{1}{r} \frac{\dd}{\dd r} ( r w_n) =0,\label{nd2}\\
& \delta^2 \frac{\dd^2  \phi}{\dd x^2}   +  \frac{1}{r} \frac{\dd}{\dd r} \left( r \frac{\dd  \phi}{\dd r} \right)  = \frac{1}{\Lambda^2} ( \ns-\ps ), \label{pois_scal} \\%\label{nd3} \\
& u_p=-\kappa_p \left( \frac{\dd p}{\dd x} + p \frac{\dd \phi}{\dd x} \right),\quad w_p=-\frac{\kappa_p}{\delta^2}\left( \frac{\dd p}{\dd r} + p \frac{\dd \phi}{\dd r} \right), \label{nd4}\\
& u_n=-\kappa_n \left( \frac{\dd n}{\dd x} -n  \frac{\dd \phi}{\dd x} \right),\quad w_n=-\frac{\kappa_n}{\delta^2}\left( \frac{\dd n}{\dd r} -n \frac{\dd \phi}{\dd r} \right), \label{nd5} \\
&\left. w_p-  \frac{d R}{d x} u_p \right|_{r=R(x)}=0, \qquad  \left. w_n- \frac{d R}{d x} u_n \right|_{r=R(x)}=0, \label{nd6} \\
&\left. \frac{\dd \phi}{\dd r} - \delta^2 \frac{d R}{d x}  \frac{\dd \phi}{\dd x}\right|_{r=R(x)}= \Upsilon \left( 1+ \delta^2 \left( \frac{d R}{d x}  \right) \right)^{1/2} \sigma(x), \label{nd7} \\
& \textrm{and the scaled boundary conditions at the ends of the pore.} \label{nd8}
 \end{align}
\end{subequations}
The dimensionless parameters in the problem are defined by
\be
\label{const}
\begin{array}{llllll}
\ds  \Upsilon= \frac{ R_0  \bar \sigma }{\varepsilon V_T } , & \ds \delta = \frac{R_0}{L}, & \ds \Lambda = \frac{L_D}{R_0},  & \ds \kappa_p=\frac{D_p}{\bar{D}}, & \ds \kappa_n=\frac{D_n}{\bar{D}}, 
\end{array}
\ee
and where $L_D$, the Debye length, is given by
\bes
\ma{L_D= \left(\frac{\varepsilon V_T } {\bar c  F } \right)^{1/2}.}
\ees
Note that $\delta$ is the aspect ratio of the pore (typically small), while $\Lambda$ measures the ratio of the Debye length of the electrolyte to the typical pore width. 
Thus where $\Lambda$ is large the pore radius is much smaller than the Debye length (which is the limit used to derive 1D Area Averaged PNP equations). However given that $L_D$ for even a 
$0.01$ Molar solution is only around $4.5$nm it is much more appropriate to consider $\Lambda=O(1)$ (or possibly even $\Lambda \ll 1$). 
The other particularly interesting parameter is $\Upsilon$; if $\Upsilon \ll 1$ then the surface charge is insufficient to induce significant ion concentration changes across 
the pore whereas if $\Upsilon=O(1)$, or greater, it  induces concentration changes that are sufficiently large to  alter the pore's macroscopic behaviour. 
Finally $\kappa_p$ and $\kappa_n$ are the dimensionless diffusivities, which assuming a sensible measure of typical diffusivity $\bar D$ is chosen will be of $O(1)$, unless the two ion diffusivities differ significantly. 
The scaled current is given by $\ds I^*=((F \bar{D} \bar{c} R_0^2)/L)I$ and can be determined from
\be
I(x,t)=\int_0^{R(x)} r (u_p-u_n) dr.  \label{nd-current}
\ee
\subsection{\ma{Parameter estimates and asymptotic limits}}\label{s:asymlim}
Nanopore devices vary in terms of length and opening radius as much as in terms of chemical composition. 
In this paper we  focus on long and narrow pores, which have been studied in many experimental setups covered in the literature, see for example \cite{pietschmann_pccp_2013,siwy_ss_2003}. 
In these pores the length is typically magnitudes of order bigger than the radius  - for example Siwy et al. \cite{siwy_ss_2003}  work with pores of 12$\mu$m length and a few nanometers radius. 
We assume that the typical length is around $L = 1\mu$m. The usual ionic
 concentration inside the pore varies from $ 0.01$  Molar up to $ 1 $ Molar.
The variations in the concentration lead to parameter  ranging from $L_D =  0.3 - 3\text{nm}$.
The opening radius may vary in the range $R_0 = 1-100$nm, hence the aspect ratio is in the range $\delta = 10^{-3} - 10^{-1}$. 
All other parameter values are listed in Table \ref{tab:scal}.
\begin{table} [h]
\begin{center}

\begin{tabular}{|l|l||l|l|} 
\hline
$K_B$ &  $1.3806504 \times 10^{-23}$ [J/K] &$\bar \sigma $ &  $ 1$  [e/nm$^2$]  = 0.16[C/m$^2$] \\
$T$ &  $300 $  [K] & $V_T $ &  $0.025$ [V]  \\
$\epsilon_0$ & $8.854187817 \times 10^{-12}$[C /(Vm)] & $\bar c $ & $  1$  [M]   \\
$\epsilon_r$ & $78.4$  &$\bar D$ &   $10^{-9}$ [m$^2$/s] \\
$\epsilon$ &  $ \epsilon_0 \epsilon_r$ & $D_p$ & $1.33$ \\
$e$ & $1.602176\times 10^{-19}$ [C] & $D_n$ & $0.79$ \\
   $\Upsilon$  &$3.4$&& \\
\hline 
\end{tabular}
\caption{Physical constants and parameters.}\label{tab:scal}
\end{center}
\end{table}
%\subsection{Asymptotic limits \label{asymlim}} 

The discussed parameter regimes motivate the following asymptotic limits. Let the dimensionless diffusivities $\kappa_p$ and $\kappa_n$ to be both $O(1)$. 
We shall only consider $\Upsilon=O(1)$, noting that the limit $\Upsilon \ll 1$ is uninteresting (because it corresponds to a  wall charge that is too small to significantly 
affect the potential and concentrations inside the pore) and that the behaviour for the limit $\Upsilon \gg 1$ can be extracted directly from the distinguished 
limit $\Upsilon=O(1)$. The size of the one remaining parameter $\Lambda$, measuring the ratio of the Debye length to the pore radius, 
determines the asymptotic structure of the solution to the PNP 
equations. In particular there are three different limits that one might wish to  consider
\begin{enumerate}[label=\roman*)]
% \begin{compactenum}[(i)]
\item $\Lambda \gg 1$, corresponding to a Debye length that it much greater than the pore 
radius, 
\item $\Lambda=O(1)$ corresponding to a Debye length that is comparable to the pore radius, and 
\item $\Lambda \ll 1$ corresponding to a Debye length much smaller than the
pore radius. 
% \end{compactenum}
\end{enumerate}
The large $\Lambda$ limit (I) has been considered in detail in a number of previous works (\eg \cite{ionreport,burger2012nonlinear}) and is only applicable 
to extremely dilute aqueous solutions and narrow pores because the Debye length $L_D$ is very small even for fairly dilute solutions (\eg 1.3nm for a 0.1M solution). 
The small $\Lambda$ limit (III) turns out to be physically rather dull because it corresponds to a Debye length that is much smaller than the pore radius meaning the the surface charge on the 
inside of the pore is effectively screened by the electrolyte and so does not significantly alter ion transport through the pore. 
Note that a similar limit was considered by Markowich in \cite{Mar} in case of the semiconductor equations with no surface charge. 
The most interesting limit, both from a mathematical and physical perspective is (II) for which $\Lambda=O(1)$. 
Furthermore we claim that this limit is a distinguished asymptotic 
limit so that the results obtained by analysing this case also provides a good description of (III) the small $\Lambda$ limits.\\

%%%%%%%%%%%%%%%%%%%%%%%%%%%%%%%%%%%%%%%%%%%%%%%%%%%%%%%%%%%%%%%%%%%%%%%

 %%%% 							Distinguished Limit

%%%%%%%%%%%%%%%%%%%%%%%%%%%%%%%%%%%%%%%%%%%%%%%%%%%%%%%%%%%%%%%%%%%%%%%

  \section{Asymptotic analysis in the  limit $\Lambda=O(1)$, $\Upsilon=O(1)$, $\delta \ll 1$ and derivation of the Quasi-1D PNP model}

  \label{medium_regime}

In this section  we discuss large aspect ratio nanopores ($\delta \ll 1$) with radii comparable  to the Debye length, (\ie  $ L_D= O(R_0)$  and hence $\Lambda =O(1)$).  In this scenario the influence of the surface charge cannot be averaged over the area of the domain, resulting in a leading order problem that must be solved both in $x$ and $r$. 
As discussed in Section \ref{s:asymlim}, we shall also consider significant surface charge by formally taking the distinguished limit $\Upsilon=O(1)$. 

%Furthermore we note that the dimensions and surface charge densities in real ion pores lead to values of $\Upsilon$, in the range 10-40. Nevertheless applications with such large values of $\Upsilon$ are still well-described by the results obtained from analysis of the distinguished limit $\Upsilon=O(1)$.

In order to find an asymptotic solution  of system (\ref{nd1})-(\ref{nd8}) in the limit $\delta \rightarrow 0$, and with all other parameters of size $O(1)$ we make the following ansatz: 
\begin{align} \label{expansion}
\begin{array}{lll}
n=n_0(r,x,t)+ \delta n_1(r,x,t)+ \cdots,  & p=p_0(r,x,t) +\delta p_1(r,x,t)+ \cdots, \\
\phi=\phi_0(r,x,t)+ \delta \phi_1(r,x,t)+ \cdots,\\
u_n= \unz(r,x,t)+ \delta \uno(r,x,t)+ \cdots, & w_n= \delta \wno(r,x,t)+ \cdots, \\
u_p= \upz(r,x,t)+\delta \upo(r,x,t)+ \cdots, & w_p= \delta \wpo(r,x,t)+ \cdots. 
\end{array} 
\end{align} 
At leading order in $\delta$ in the flux equations (\ref{nd4})-(\ref{nd5}) give the two equations
\begin{align*}
 \frac{\partial n_0}{\partial r}  - n_0  \frac{\partial \phi_0}{\partial r}= 0 ~~ \text{ and }  ~~ \frac{\partial p_0}{\partial r}  +p_0  \frac{\partial \phi_0}{\partial r}= 0 ,
    \end{align*}
which can be integrated to obtain 
\be
n_0 = Q(x,t) \exp(\phi_0(r,x,t)) ~~\text{ and }  ~~p_0 = S(x,t) \exp(-\phi_0(r,x,t)), \label{npsol}
\ee
where the functions $Q(x,t)$ and $S(x,t)$ are yet to be determined.  
Inserting these expressions into the Poisson-Boltzmann equation and its boundary condition, (\ref{pois_scal}) and (\ref{nd7}) gives 
\begin{subequations}\label{poiss_1D}
\begin{align}
\frac{1}{r} \frac{\dd}{\dd r} \left( r \frac{\partial \phi_0}{\partial r} \right) &=  \frac{1}{\Lambda^2} \left(Q(x,t) \exp(\phi_0(r,x,t)) -  S(x,t) \exp(-\phi_0(r,x,t)) \right),  \label{poteq1}\\
 \left. \frac{\partial   \phi_0}{\partial r  }  \right|_{r=R(x)}&=   \Upsilon \sigma(x).~~~~~~~~~~~~~~~~~~~~~~~~~~~~~~~~~~~  \label{poteq2}
 \end{align}
\end{subequations}

 \subsection{Leading order solution for the potential \label{lo-pot}}
We seek solutions to (\ref{poteq1})-(\ref{poteq2}) by introducing the new variables (\eg see \cite[Chapter 12]{andelman1995electrostatic})
\be
\phi_0(r,x,t)=  \frac{1}{2} \log{\frac{S(x,t)}{Q(x,t)}}  + \psi (\re{\xi},x,t) \quad   \text{and } \quad r = R(x) \xi, \label{kim}
\ee 
which result in the following problem for $\psi$:
\begin{subequations}\label{sys_sinh}
\begin{align}
&\frac{1}{\xi} \frac{\dd}{\dd \xi} \left( \xi \frac{\partial \psi}{\partial \xi} \right) = \frac{1}{\lambda^2(x,t)} (e^{\psi}-e^{-\psi}),  \label{sinh} \\ 
&\psi \ \ \mbox{bounded} \ \mbox{at} \ \ \xi=0, \textrm{ and } \frac{\partial \psi   }{\partial \xi  } \Big  |_{\xi =1}  =  \beta(x).  \label{sinh_bc}
\end{align}
\end{subequations}
where
\be
\beta(x)=\Upsilon \sigma(x) R(x), \qquad \mbox{and} \qquad \lambda(x,t)= \frac{\Lambda}{R(x) (S(x,t)Q(x,t))^{1/4}}. \label{lambet}
\ee
Here $\lambda(x,t)$ gives the ratio of the Debye length, evaluated from the evolving ion concentrations, to the local pore radius. 
Thus the solution to  (\ref{sinh})-(\ref{sinh_bc}),
\bes
\psi=\psi \left(\xi; \lambda(x,t),  \beta(x) \right),
\ees
 depends parametrically on $x$ and $t$ through $\lambda(x,t)$ and $\beta(x)$.
 
 \paragraph{{Approximate solution to  (\ref{sinh})-(\ref{sinh_bc}) for $\beta(x) \gg 1$.}}
We can make use of the fact that $\beta(x)$ is typically large (so that the gradient of $\psi$ at the edge of the pore $\xi=1$ is large) by noting that this suggests that $\psi(\xi; \lambda,  \beta)$ is also large (an hypothesis we justify {\it a posteriori} for sufficiently large $\lambda$). Making the large $\psi$ ansatz means that (\ref{sinh}) can be approximated by
\be
\frac{1}{\xi} \frac{\dd}{\dd \xi} \left( \xi \frac{\partial \psi}{\partial \xi} \right) \sim \frac{1}{\lambda^2(x,t)} e^{\psi} \label{sinh_approx}
\ee
which, when solved together with (\ref{sinh_bc}), has a  solution of the form
\be
\psi(\xi,x,t) \sim 2 \log \left( \mbox{cosech} \left( \mbox{arcoth} \left( \frac{\beta(x)+2}{2} \right) -\log \xi \right) \right) -
\log \left( \frac{\xi^2}{2 \lambda^2(x,t)} \right).  \label{psi_approx}
\ee
Notably this expression for $\psi$ has a minimum (as a function of $\xi$) at the centre of the pore given by
\be
\psi|_{\xi=0}=3 \log 2 + 2 \log \lambda- 2 \mbox{arcoth} \left( \frac{\beta+2}{2} \right) ,
\ee
which for $\beta \gg 1$ is well-approximated by $\psi|_{\xi=0}=3 \log 2 + 2 \log \lambda$. The approximation in going from (\ref{sinh}) to (\ref{sinh_approx}) can thus be justified if $\exp(-2 \psi|_{\xi=0}) \ll 1$ which is true only if
\bes
\lambda(x,t) \gg \frac{1}{2^{3/2}}.
\ees

\begin{figure}[H]
    \centering
    \begin{subfigure}[b]{0.485\textwidth}
        \includegraphics[width=\textwidth]{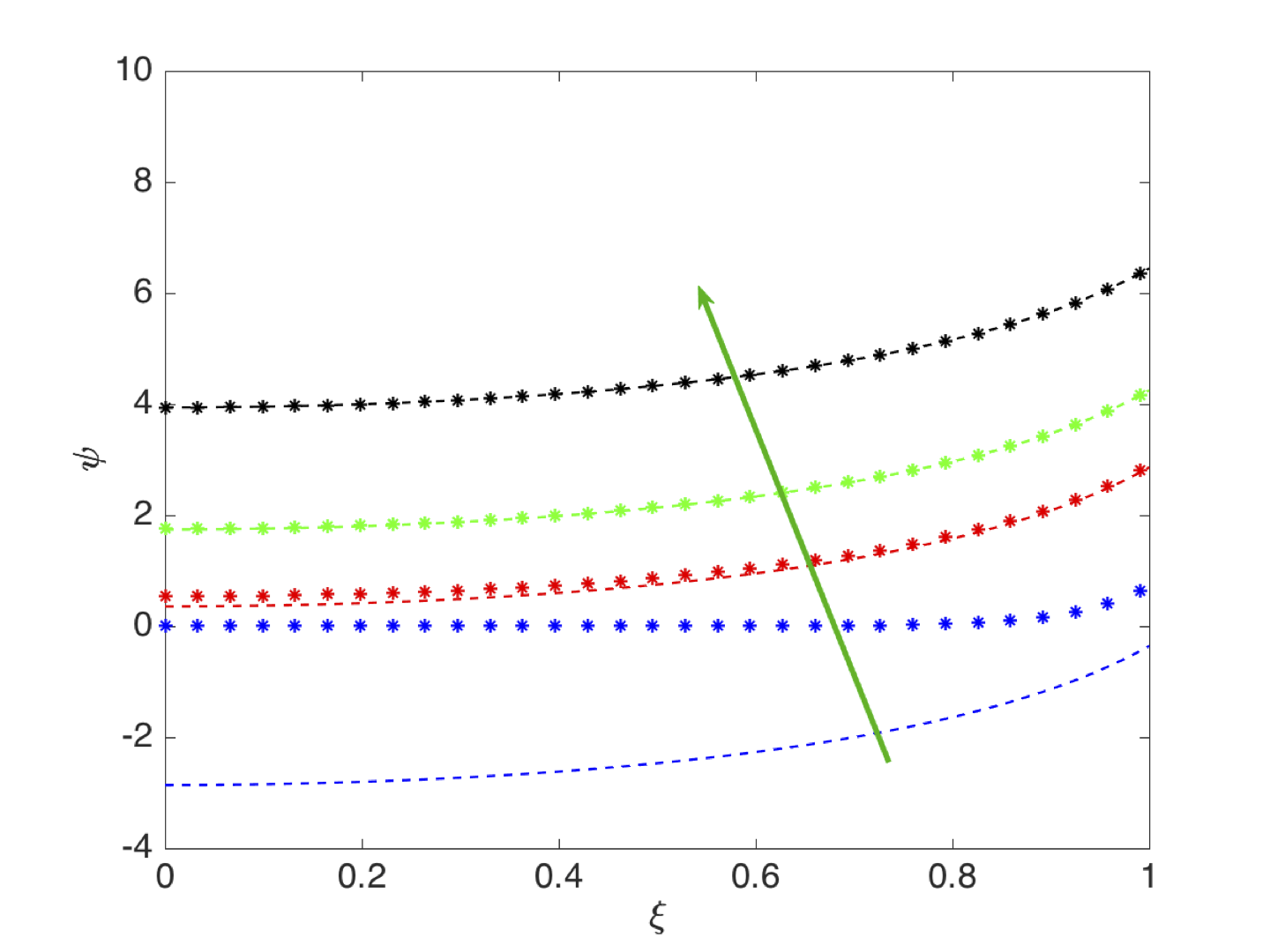}
    \end{subfigure}
    ~ 
    \begin{subfigure}[b]{0.485\textwidth}
        \includegraphics[width=\textwidth]{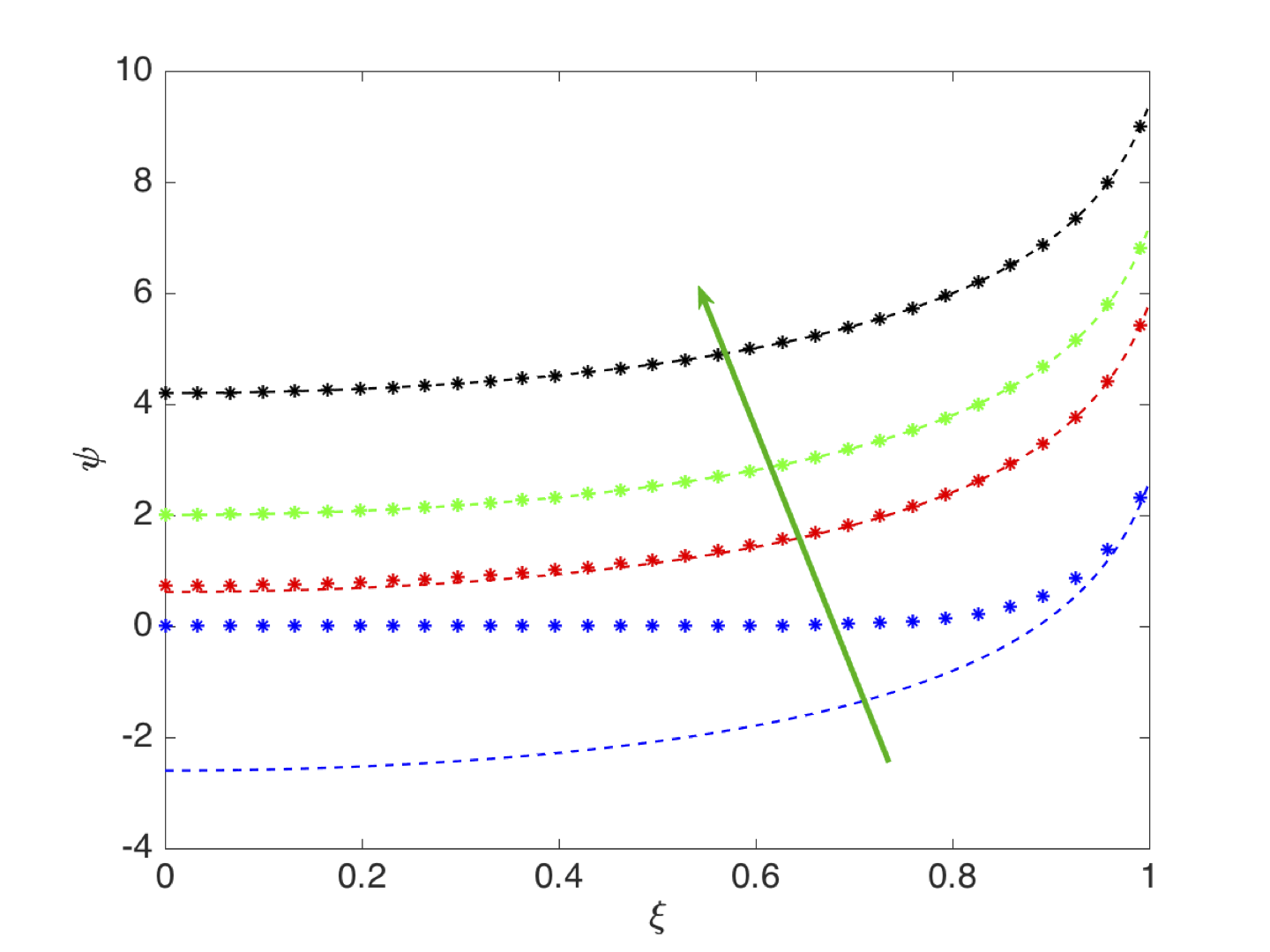}
      \end{subfigure}
       \caption{Comparison between numerical solution to \eqref{sys_sinh} (stars) and its large-$\beta$ asymptotic approximation \eqref{psi_approx}  (dashed line). Here in the left-hand panel  $\beta=10$ while in the right $\beta=50$. In both panels the values of $\lambda$ taken are $\lambda = [0.1,0.5,1,3]$ and the arrows indicate the direction of increasing $\lambda$.}
       \label{bigbeta_comparisson}
\end{figure}

Figure \ref{bigbeta_comparisson} illustrates the quality of the asymptotic solution for different values of $\lambda$ and realistic values of $\beta=10,50$.  
Note that the approximation quality of the asymptotic solution improves as $\lambda$ increases.

\paragraph{Approximate solution to  (\ref{sinh})-(\ref{sinh_bc}) for $\lambda(x,t) \ll 1$.}
In this regime narrow Debye layers of width $O(\lambda)$ exist close to surface $\xi=1$. In order to investigate the solution further we  rescale, in the standard manner (see \cite{markowich1986asymptotic,markowich1986uniform}), about this surface by introducing the Debye layer coordinate $\z$ defined by 
\be
\xi=1-\lambda \z.
\ee
Rewriting (\ref{sinh})-(\ref{sinh_bc}) in terms of this new coordinate leads to the following equation and boundary condition for $\psi$
\be
\frac{\dd^2 \psi}{\dd \z^2} + \lambda \left( \z \frac{\dd^2 \psi}{\dd \z^2} - \frac{\dd \psi}{\dd \z} \right) + O(\lambda^2)= (e^{\psi}-e^{-\psi}), \label{chorizo1}\\
\left. \frac{\dd \psi}{\dd \z} \right|_{\z=0} = -B \qquad \mbox{where} \quad B=\lambda \beta,\label{chorizo2}\\
\frac{\dd \psi}{\dd \z} \ra 0 \qquad \mbox{as} \quad \z \ra +\infty. \label{chorizo3}
\ee
Here we consider the distinguished limit $B=O(1)$, that is $\beta=O(1/\lambda)$ noting that the solution we obtain is still valid for other sizes of this parameter. Formally we look for a solution in the Debye layer by expanding $\psi$ in the form $\psi=\psdz+\lambda \psdo +\cdots$, substituting into 
(\ref{chorizo1})-(\ref{chorizo3}) and taking the leading order terms. This results in the following problem for $\psdz$
\be
\frac{\dd^2 \psdz}{\dd \z^2} =(e^{\psdz}-e^{-\psdz}),~~~~~~~~~~~~~~~~~~~~~~~~~~ \\
\left. \frac{\dd \psdz}{\dd \z} \right|_{\z=0} = -B, \qquad \mbox{and} \qquad \frac{\dd \psdz}{\dd \z} \ra 0 \quad \mbox{as} \quad \z \ra +\infty.
\ee
This, as is well-known, has the solution
\be
\psdz=\left\{ \begin{array}{lcc} \ds  2 \log_e \left( \coth \left[ \frac{1}{\sqrt{2}} \left( \z + \frac{1}{\sqrt{2}} \mbox{arcsinh} \left( \frac{2 \sqrt{2}}{B} \right) \right)\right] \right) & \mbox{for} & B>0, \\*[6mm]
\ds  2 \log_e \left( \tanh \left[ \frac{1}{\sqrt{2}} \left( \z + \frac{1}{\sqrt{2}} \mbox{arcsinh} \left( -\frac{2 \sqrt{2}}{B} \right) \right)\right] \right) & \mbox{for} & B<0.
\end{array} \right.  \label{custard}
\ee
Notably this solution has the property that
\bes
\psdz \ra 0 \quad \mbox{as} \quad \z \ra +\infty,
\ees
and so is uniformly valid for all values of $\xi \in [0,1)$ or equivalently for $\z \in [1/\lambda,0)$. It follows that we do not need to look for a solution for $\psi$ in an outer region.
 
\begin{figure}[H]
    \centering
    \begin{subfigure}[b]{0.485\textwidth}
        \includegraphics[width=\textwidth]{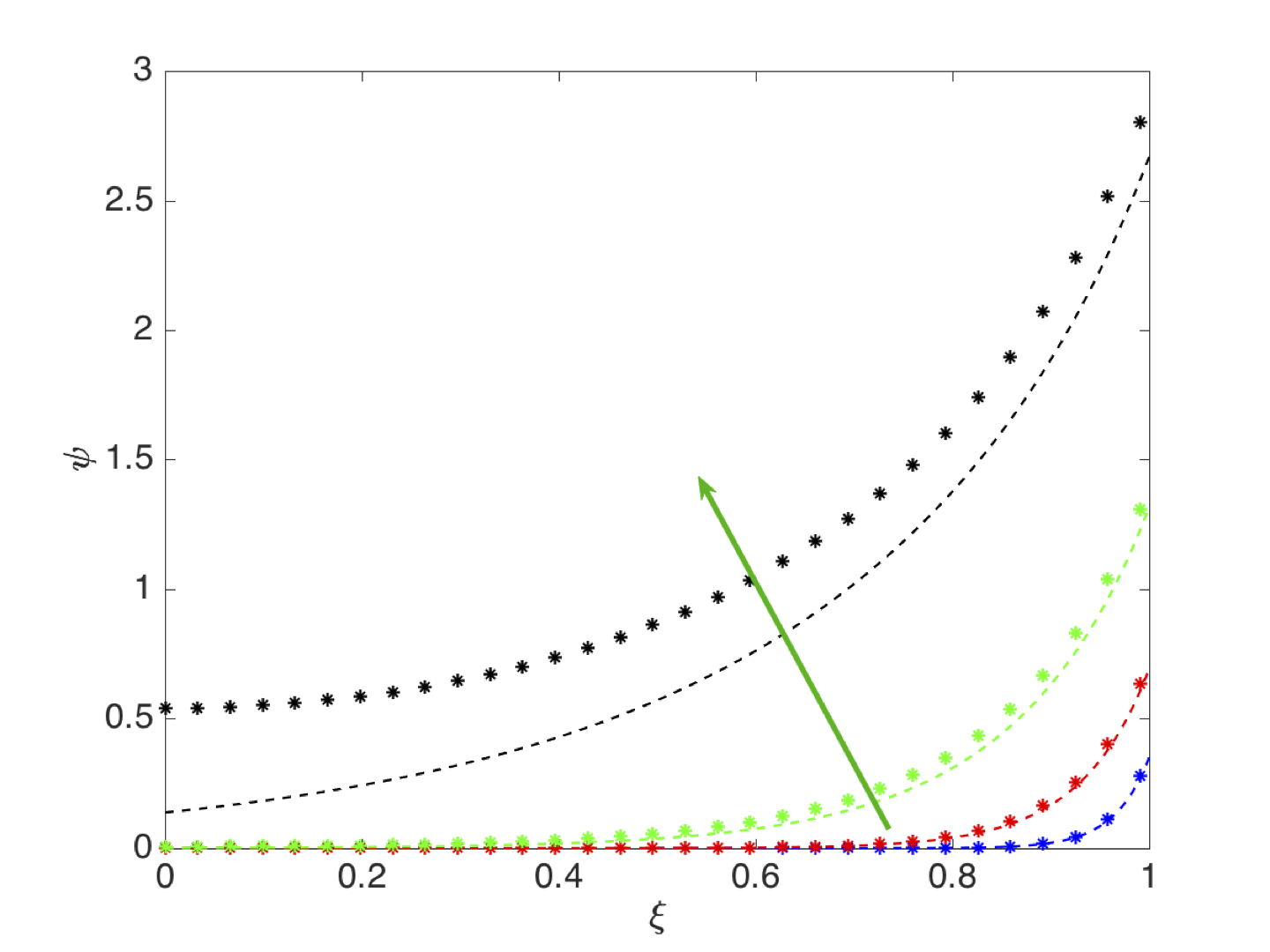}
    \end{subfigure}
    ~ 
    \begin{subfigure}[b]{0.485\textwidth}
        \includegraphics[width=\textwidth]{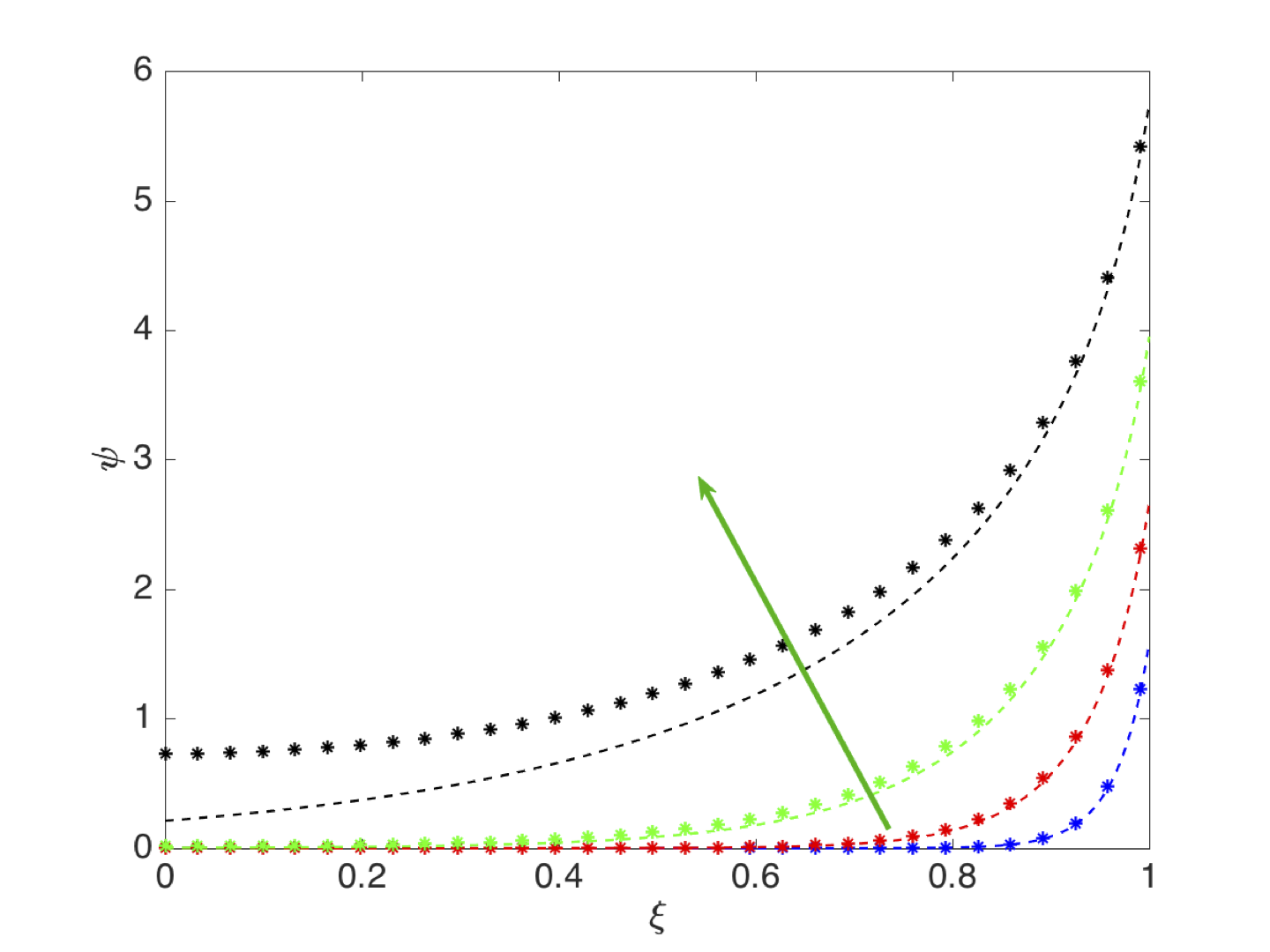}
      \end{subfigure}
       \caption{Comparison between numerical solution to \eqref{sys_sinh} (stars) and its small-$\lambda$ asymptotic approximation \eqref{psi_approx}  (dashed line). Here in the left-hand panel  $\beta=10$ while in the right $\beta=50$. In both panels the values of $\lambda$ taken are  $\lambda = [ 0.05, 0.1 ,0.2, 0.5]$ and the arrows indicate the direction of increasing $\lambda$.}
       \label{small_lambda}
\end{figure} 
Figure \ref{small_lambda} illustrates the asymptotic solution as well as the numerical solution of \eqref{sys_sinh} for different values of $\lambda$ and realistic values of $\beta=10,50$.  
In this case the approximation quality, as expected, increases as $\lambda$ decreases.
%  
%  It remains to find equations for the evolution of $Q(x,t)$ and $S(x,t)$.
% 
%%%%%%%%%%%%%%%%%%%%%%%%%%%%%%%%%%%%%%%%%%%%%%%%%%%%%%%%%%%%%%%%%%%%%
%%%%%  FLUX CONSERVATION AND ONE-D MODEL %%%%%%%%%%%%%%%%%%%%%%%%%%%%%%%%%%%%%%%%
%%%%%%%%%%%%%%%%%%%%%%%%%%%%%%%%%%%%%%%%%%%%%%%%%%%%%%%%%%%%%%%%%%%%%

 \subsection{Leading order flux conservation and the simplified 1D model \label{lo-flux}}
The purpose of this section is to derive flux conservation conditions in the $x$-direction that will give rise to evolution equations for $Q(x,t)$ and $S(x,t)$. 
Together with the approximations of the previous section, we are then in a position  to numerically solve the approximated system. We start by considering  the leading order terms in the ion conservation equations (\ref{nd1})-(\ref{nd2}), namely
\be
\frac{\dd p_0}{\dd t}+\frac{\partial u_{p,0}}{\partial x} + \frac{1}{r} \frac{\partial }{\partial r } (r w_{p,1})=0, \label{huevo1} \\
\frac{\dd n_0}{\dd t}+\frac{\partial u_{n,0}}{\partial x} + \frac{1}{r} \frac{\partial }{\partial r } (r w_{n,1})=0. \label{huevo2}
\ee
where expressions for $u_{p,0}$ and $u_{n,0}$ are obtained from the leading order expansions of (\ref{nd4}) and (\ref{nd5}) and are
\be
u_{p,0}= -\kappa_p \left( \frac{\dd p_0}{\dd x} + p_0 \frac{\dd \phi_0}{\dd x}  \right), \quad \mbox{and} \quad u_{n,0}= -\kappa_n \left( \frac{\dd n_0}{\dd x} -n_0 \frac{\dd \phi_0}{\dd x}  \right).  \label{unamo}
\ee
The boundary conditions on $w_{p,1}$  and $w_{n,1}$  come from the leading order expansion of (\ref{nd6}) and are
\be
w_{p,1} |_{r=R(x)}= \left. \frac{d R}{d x}  u_{p,0} \right|_{r=R(x)}, \qquad w_{n,1} |_{r=R(x)}= \left. \frac{d R}{d x}  u_{n,0} \right|_{r=R(x)}. \label{huevo3}
\ee
Multiplying both (\ref{huevo1})  and (\ref{huevo2}) by $r$ and integrating between $r=0$ and $r=R(x)$ gives
\begin{align*}
\int_0^{R(x)} \left(\frac{\dd p_0}{\dd t}+\frac{\partial u_{p,0}}{\partial x}  \right) r dr + \left[ r w_{p,1} \right]_{r=0}^{R(x)} =0, \quad \int_0^{R(x)} \left(\frac{\dd n_0}{\dd t}+\frac{\partial u_{n,0}}{\partial x}  \right) r dr + \left[ r w_{n,1} \right]_{r=0}^{R(x)} =0.
\end{align*}
On applying the boundary conditions (\ref{huevo3}) it can be seen that these equations can be rewritten in conservation form
\be
\frac{\dd}{\dd t} \left(\int_0^{R(x)}r  p_0(r,x,t) dr \right) + \frac{\partial }{\partial x} \left(\int_0^{R(x)} r u_{p,0} dr  \right)  &=&0,\\
\frac{\dd}{\dd t} \left(\int_0^{R(x)}r  n_0(r,x,t) dr \right) + \frac{\partial }{\partial x} \left(\int_0^{R(x)} r u_{n,0} dr  \right)  &=&0.
\ee
Substituting for $p_0(r,x,t)$ and $n_0(r,x,t)$ from (\ref{npsol}) and $u_{p,0}$ and $u_{n,0}$ from (\ref{unamo}) leads to an alternative reformulation
\be
\frac{\dd}{\dd t} \left(S(x,t) \Theta_1(x,t) \right) &=& \kappa_p  \frac{\partial }{\partial x} \left( \frac{\dd S}{\dd x} \Theta_1(x,t) \right),  \label{tory1}\\
\frac{\dd}{\dd t} \left(Q(x,t) \Theta_2(x,t) \right) &=& \kappa_n  \frac{\partial }{\partial x} \left( \frac{\dd Q}{\dd x} \Theta_2(x,t) \right), \label{tory2}
\ee
where
\be
\Theta_1(x,t)=\pi \int_{r=0}^{R(x)} r \exp( -\phi_0(r,x,t)) dr,\quad \mbox{and} \quad \Theta_2(x,t)=\pi \int_{r=0}^{R(x)} r \exp( \phi_0(r,x,t)) dr. \label{Thetadef}
\ee
On substituting for $\phi_0$ and $r$, in terms of $\psi$ and $\xi$, from (\ref{kim}) we can rewrite these expressions in the form
\be
\Theta_1(x,t)=A(x) \left(\frac{Q(x,t)}{S(x,t)}\right)^{1/2} G_1(x,t), \qquad  \Theta_2(x,t)=A(x)\left(\frac{S(x,t)}{Q(x,t)}\right)^{1/2} G_2(x,t),  \label{Theta}
\ee
where $A(x)=\pi R^2(x)$ is the cross-sectional area of the pore and the functions $G_1$ and $G_2$ are defined by
\be
G_1(x,t)&=&\int_{\xi=0}^{1} \xi \exp \left( -\psi \left(\xi; \lambda(x,t),  \beta(x) \right) \right) d \xi,  \label{G1}\\
 G_2(x,t)&=&  \int_{\xi=0}^{1} \xi \exp \left( \psi \left(\xi; \lambda(x,t),  \beta(x) \right) \right)  d \xi.~~~~~\label{G2}
\ee
Here $\lambda(x,t)$ and $\beta(x)$ are defined in (\ref{lambet}). Notably since $\psi(x,t)$ satisfies the problem (\ref{sinh})-(\ref{sinh_bc}) we can show (by multiplying (\ref{sinh}) by $\xi \lambda^2(x,t)$, integrating between $\xi=0$ and 1 and imposing the boundary conditions (\ref{sinh_bc})) that
\be
G_2(x,t)-G_1(x,t)=\lambda^2(x,t) \beta(x). \label{Gdiff}
\ee

The leading order current flowing through the pore can be calculated from (\ref{nd-current}), (\ref{npsol}), (\ref{unamo}) and (\ref{Thetadef}) and is
\be
I \sim -\kappa_p \frac{\dd S}{\dd x}  \Theta_1+\kappa_n \frac{\dd Q}{\dd x}  \Theta_2. \label{leadflux}
\ee
or equivalently, on referring to (\ref{Theta}),
\be
I \sim A(x) ( S Q)^{1/2} \left( \kappa_n G_2 \frac{\dd}{\dd x} \log_e Q - \kappa_p G_1 \frac{\dd}{\dd x} \log_e S \right). \label{currSQ} 
\ee
\begin{rem} \label{rem:decoupled}
As mentioned in the introduction, the most commonly available data from nanopore experiments are IV curves. 
Thus equation \eqref{leadflux} (and \eqref{currSQ}) allow us to compute the IV curves very efficiently, without solving a non-linear Poisson equation 
(as it is the case in the classical Scharfetter--Gummel iteration for PNP, \cite{gummel1964self}). From the computational point of view this is the main advantage of our approach.
\end{rem}

\begin{figure}
     \begin{center}
       \includegraphics[width=1\textwidth]{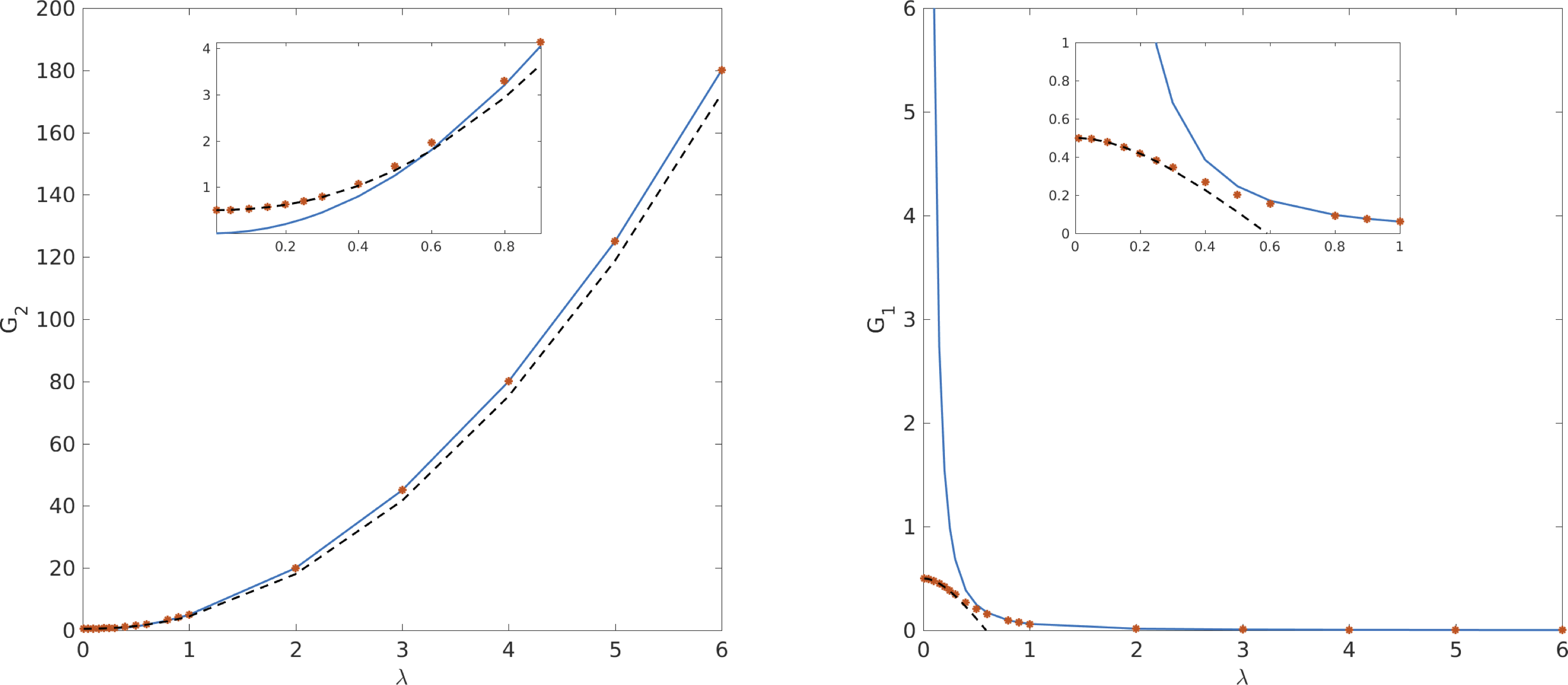}
        \end{center}
\begin{center}
        \includegraphics[width=1\textwidth]{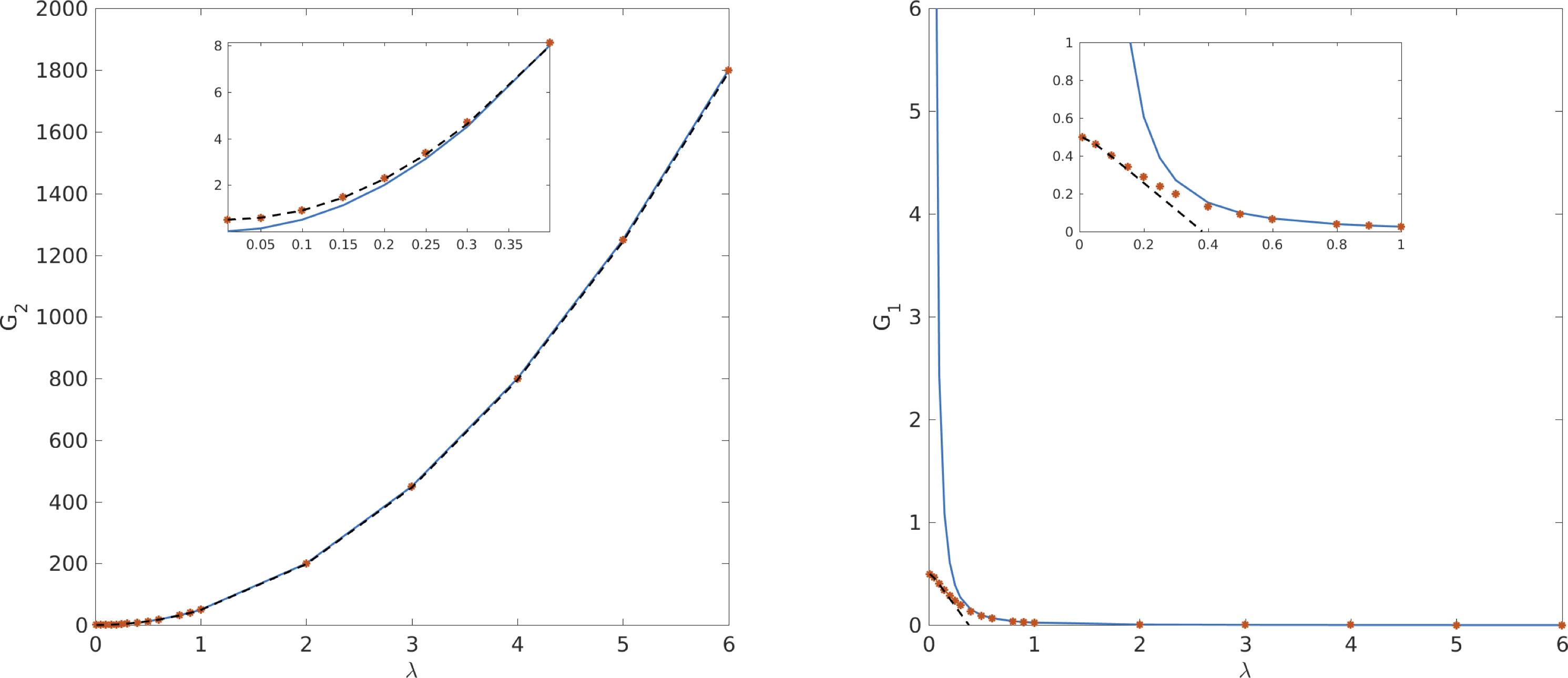}
               \end{center}
\caption{Comparison of the numerical evaluation of the expressions $G_1$  (right) and $G_2$  (left) as a function of $\lambda$ for $\beta = 5$ (top row) and $\beta = 50$ (bottom row). Red stars correspond to the values of $G_1$ and $G_2$ calculated from the full equations \eqref{G1}--\eqref{G2} The blue solid line corresponds to the approximation \eqref{G1G2bigbeta} ($\beta \gg 1$) while the black dashed one stands for \eqref{G1G2smalllambda} ($\lambda \ll 1$).}                 
               \label{G1G2comp}
\end{figure}

\paragraph{Approximations of $G_1(x,t)$ and $G_2(x,t)$ for $\beta \gg 1$ and $\lambda=O(1)$.} 
 In this instance we can find asymptotic expressions for $G_1$ and $G_2$ simply by substituting (\ref{psi_approx}), the large $\beta$ asymptotic expression for $\psi$, directly into (\ref{G1})-(\ref{G2}) to obtain 
\be
G_1(x,t) \sim \frac{1}{48 \lambda^2(x,t)} \frac{\beta^2(x)+12 \beta(x) + 48}{\beta(x) (\beta(x)+4)}, \quad \mbox{and} \quad G_2(x,t) \sim \lambda^2(x,t) \beta(x). 
\label{G1G2bigbeta}
\ee
Note that these expressions still satisfy the identity (\ref{Gdiff}) asymptotically in the limit $\beta \ra \infty$ since $G_1 \ll G_2$. However the asymptotic expansion breaks down for $\lambda \ll 1$, as noted previously, and so we need to obtain alternative expressions for $G_1$ and $G_2$ in this limit.
\begin{rem}
Note that the case   $\beta \ll -1$  can be solved by setting  $ u = -\psi$  in equation \eqref{sinh_approx} and following the calculations detailed above to obtain approximations for $ G_1$ and $G_2$. 
\end{rem}
\paragraph{Approximations of $G_1(x,t)$ and $G_2(x,t)$ for $\lambda \ll 1$.}
In order to approximate the integrals in (\ref{G1})-(\ref{G2}) based on the Debye layer solution for $\psi$ in the small $\lambda$ limit (\ref{custard}) we split the integrals up as follows
\begin{align*}
 G_1(x,t)&=\int_{\xi=0}^{1}\xi d\xi- \int_{\xi=0}^{1} \xi \left(1-\exp \left( -\psi \left(\xi; \lambda,  \beta \right) \right) \right)d \xi, \\
 G_2(x,t)&=\int_{\xi=0}^{1}\xi d\xi- \int_{\xi=0}^{1} \xi \left(\exp \left( \psi \left(\xi; \lambda,  \beta \right) \right)-1 \right)  d \xi
\end{align*}
before substituting $\xi=1-\lambda \z$ and formally taking the limit $\lambda \ra 0$ to obtain the following asymptotic expressions
\bes
G_1 \sim\frac{1}{2}- \lambda \int_{\z=0}^{\infty}  \left(1-\exp \left( -\psdz \right) \right)d \z\quad \mbox{and} \quad G_2 \sim \frac{1}{2}+ \lambda \int_{\z=0}^{\infty}  \left(\exp \left( \psdz \right) -1 \right)d \z.
\ees
Evaluating these expressions, in the distinguished limit that $B=\lambda \beta=O(1)$, gives the following relations for $G_1$ and $G_2$ in the small $\lambda$ limit
\be
G_1 \sim \frac{1}{2} -\lambda \frac{2 \sqrt{2} B}{\sqrt{8+B^2}+2 \sqrt{2} +B}\quad \mbox{and} \quad G_2\sim \frac{1}{2} +\lambda \frac{2 \sqrt{2} B}{\sqrt{8+B^2}+2 \sqrt{2} -B}.
\label{G1G2smalllambda}
\ee
In this instance it turns out that these asymptotic expressions for $G_1$ and $G_2$, which are formally of the same order, satisfy the condition (\ref{Gdiff}) identically.
Figure \ref{G1G2comp} shows that by choosing the right cut-off value  of $\lambda$  it is possible to obtain an adequate approximations to $G_1$ and $G_2$ for all values of $\lambda$ provided that $\beta$ is large.  This approximation can be much improved by smoothing between the two asymptotic representations of the solutions, in the limits $\beta \gg1$ and $\lambda \ll 1$.  The smoothed, uniformly valid asymptotic, representation of $G_1$ is discussed further in Appendix \ref{B} and the accuracy of the fit to the numerical solutions for $G_1$ can be appreciated by  inspecting Figure \ref{G1_smooth}. Note that once a good representation of $G_1$ has been obtained $G_2$ can be directly evaluated from the relation (\ref{Gdiff}). \\

\subsection{Summary of the Quasi-1D model} 
 Since the resulting 1D model (comprised of equations (\ref{npsol}), (\ref{kim})-(\ref{lambet}), (\ref{tory1})-(\ref{tory2}), (\ref{Theta})-(\ref{G1}), 
 (\ref{Gdiff}) and \eqref{G1G2smalllambda}, \eqref{G1G2bigbeta}  is quite intricate we summarise it in the following paragraph. 
 The leading order ion concentrations and potential are given in terms of the functions $Q(x,t)$, $S(x,t)$ and $\psi(\xi,x,t)$ by the following:
\be
n(r,x,t)&=&(Q(x,t) S(x,t))^{1/2} \exp\left(\psi \left(\frac{r}{R(x)},x,t \right) \right), \label{ne} \\ 
p(r,x,t)&=&(Q(x,t) S(x,t))^{1/2} \exp\left(-\psi \left(\frac{r}{R(x)},x,t \right) \right), \label{pe} \\
\phi(r,x,t)&=&\frac{1}{2} \frac{S(x,t)}{Q(x,t)} +\psi \left(\frac{r}{R(x)},x,t \right), \label{phie} 
\ee
where $\psi(\xi,x,t)$ 
%(note the change of independent variables from $(\xi;\lambda(x,t),\beta(x))$ to $(\xi,x,t)$) 
satisfies the following series of ODE problems in $\xi$
\be
\frac{1}{\xi} \frac{\dd}{\dd \xi} \left( \xi \frac{\partial \psi}{\partial \xi} \right) = \frac{R^2(x) (Q(x,t) S(x,t))^{1/2}}{\Lambda^2} (e^{\psi}-e^{-\psi}), \label{psi1} \\ 
\psi \ \ \mbox{bounded} \ \mbox{at} \ \ \xi=0, \textrm{ and } \frac{\partial \psi   }{\partial \xi  } \Big  |_{\xi =1}  =  \Upsilon \sigma(x) R(x).   \label{psi2}
\ee
In turn the functions $Q(x,t)$ and $S(x,t)$ satisfy the PDEs
\be
\frac{\dd}{\dd t} \left( (Q(x,t) S(x,t))^{1/2} G_1(x,t) \right) = \frac{\kappa_p}{A(x)} \frac{\dd}{\dd x} \left( A(x) 
\left(\frac{Q(x,t)}{S(x,t)} \right)^{1/2} G_1(x,t) \frac{\dd S}{\dd x} \right), \label{QQQ}\\
\frac{\dd}{\dd t} \left( (Q(x,t) S(x,t))^{1/2} G_2(x,t) \right) = \frac{\kappa_n}{A(x)} \frac{\dd}{\dd x} \left( A(x) 
\left(\frac{S(x,t)}{Q(x,t)} \right)^{1/2} G_2(x,t) \frac{\dd Q}{\dd x} \right), \label{SSS}
\ee
where $G_2(x,t)$ and $G_1(x,t)$ are given by the expressions
\be
G_2(x,t) &=& G_1(x,t) +\frac{\Lambda^2 \Upsilon \sigma(x)}{R(x) (Q(x,t) S(x,t))^{1/2}}, \label{GG2} \\
G_1(x,t) &=& \int_0^1 \xi \exp( -\psi(\xi,x,t) ) d \xi.  \label{GG1}
\ee
The main point of the method is that the integrals $G_1$ and $G_2$ are not calculated via integrating the  $\psi$ but using the polynomial approximations obtained 
in the equations \eqref{G1G2smalllambda} and \eqref{G1G2bigbeta} for different values  of the $\lambda(x)$. Thus \eqref{QQQ}--\eqref{SSS} are decoupled 
from \eqref{psi1}--\eqref{psi2}. As mentioned in Remark \ref{rem:decoupled}, this is a particular advantage since it allows to calculate the ion current 
(via \eqref{currSQ}) without having to solve a nonlinear equation.
To ensure a smooth transition between the two regimes a smoothing procedure was implemented   as described in detail in the  Appendix \ref{B}, that is by writing
\be
G_1(x,t) \approx \gsm \left(\frac{\Lambda}{R(x) (S(x,t)Q(x,t))^{1/4}}, \Upsilon \sigma(x) R(x) \right)
\ee
where the function $\gsm(\lambda,\beta)$ is defined in (\ref{gsm}). This approximation of $G_1(x,t)$ taken together with (\ref{QQQ})-(\ref{GG2}) allows us
to solve a one-dimensional spatial problem for $S(x,t)$ and $Q(x,t)$. 
If the purpose of the calculation is solely to determine the current $I$  flowing through the pore (for example when calculating $I$-$V$ curves) this calculation is 
sufficient since $I$ may be calculated solely from $S(x,t)$, $Q(x,t)$, $G_1(x,t)$ and $G_2(x,t)$ via the formula \eqref{currSQ}, that is by
\be
I = A(x) ( S Q)^{1/2} \left( \kappa_n G_2 \frac{\dd}{\dd x} \log_e Q - \kappa_p G_1 \frac{\dd}{\dd x} \log_e S \right). 
\ee
If in addition to determining the current flow through the pore we wish also to obtain the spatial distributions of the carrier concentrations and the electric potential  we need also to solve for the function $\psi(\xi,x,t)$ in order to use it in \eqref{ne}-\eqref{phie} in order to calculate
$n(r,x,t)$, $p(r,x,t)$ and $\phi(r,x,t)$. Although it is possible to obtain a reasonable approximation to the function $\psi(\xi,x,t)$ in the large $\beta$ limit by using the appropriate asymptotic solution, \eqref{psi_approx} for $\lambda=O(1)$ or \eqref{custard} for $\lambda \ll 1$, we instead choose to  solve  the full boundary value problem for $\psi(\xi,x,t)$ numerically, as specified in \eqref{sinh}-\eqref{sinh_bc}; that is we solve
\be
\frac{1}{\xi} \frac{\dd}{\dd \xi} \left( \xi \frac{\partial \psi}{\partial \xi} \right) = \frac{R^2(x) (S(x,t)Q(x,t))^{1/2}}{\Lambda^2}  (e^{\psi}-e^{-\psi}),  \label{sinh2} \\ 
\psi \ \ \mbox{bounded} \ \mbox{at} \ \ \xi=0, \textrm{ and } \left. \frac{\partial \psi   }{\partial \xi  } \right|_{\xi =1}  =  \Upsilon \sigma(x) R(x),  \label{sinh_bc2}
\ee
for each position $x$ and time $t$. We adopt this numerically costly procedure here because it provides more accurate asymptotic representations of 
$n(r,x,t)$, $p(r,x,t)$ and $\phi(r,x,t)$ with which to compare to the full 2D numerical solutions 
(see figures \ref{heatmap_asym}, \ref{profiles}, \ref{heatmap_bigpore_2D_asym} and \ref{profiles_bigpore}). 
We do however believe that it should be possible to obtain a uniformly valid asymptotic expansion 
for $\psi(\xi,x,t)$ in the large $\beta$ limit that is capable of accurately capturing the solution for all values of $\lambda$, much as we do for $G_1$ in Appendix \ref{B}.

\subsection{An alternative formulation.}
It is possible to reformulate the Quasi One-1D PNP Model derived in \S \ref{lo-pot}-\ref{lo-flux} and contained in (\ref{tory1})-(\ref{G2}) in more physically appealing forms. We give one such reformulation below but note that there are others.

We start by noting that the (dimensionless) electrochemical potentials of positive and negative ions, $\mu_p$ and $\mu_n$ respectively, are 
\be
\mu_p=\log_e p+ \phi, \quad \mbox{and} \quad \mu_n=\log_e  n - \phi.
\ee
The chemical potential of the electrolyte, defined by $\mu_e=\frac{1}{2} ( \mu_p+\mu_n)=\frac{1}{2} \log_e(np)$, is obtained at leading order by substituting the approximations to $n$ and $p$ found in (\ref{npsol}) into this expression; this gives
\be
\mu_e(x,t)=\log \left(\left(Q(x,t) S(x,t) \right)^{1/2} \right).  \label{chempot}
\ee
In addition we define an effective electric potential $\phit$ by
\be
\phit(x,t)= \log_e \left( \left( \frac{S(x,t)}{Q(x,t)} \right)^{1/2} \right). \label{effpot}
\ee
We now introduce two further quantities $\pb$ and $\nb$, the cross-sectionally averaged ion densities, as defined by
\be
\pb(x,t)=\frac{\pi}{A(x)} \int_0^{R(x)} r p_0 dr \quad \mbox{and} \quad \nb(x,t)=\frac{\pi}{A(x)} \int_0^{R(x)} r n_0 dr.
\ee
Substituting for $n_0$ and $p_0$ from (\ref{npsol}), and making use of the definitions (\ref{Thetadef}), allows us to re-express these quantities in the form
\be
\pb(x,t)=\frac{S(x,t) \Theta_1(x,t)}{A(x)}, \quad \mbox{and} \quad \nb(x,t)=\frac{Q(x,t) \Theta_2(x,t)}{A(x)}.  \label{NPexpr5}
\ee
In turn substituting for $\Theta_1$ and $\Theta_2$ from (\ref{Theta}), and using the formula (\ref{chempot}) to eliminate $Q$ and $S$, allows us to rewrite  $\pb$ and $\nb$ as follows:
\be
\pb(x,t)=\exp(\mu_e) G_1(\mu_e; x), \quad  \nb(x,t)= \exp(\mu_e) G_2(\mu_e;x).~~~~  \label{NPexpr4}
\ee
In the above we  have written both $G_1$ and $G_2$ in a form that makes it explicit that these quantities are independent  of $\phit$ and depend only on $Q$ and $S$ through $\mu_e(x,t)$. On using the definitions (\ref{chempot}) and (\ref{effpot}) to eliminate $Q$ and $S$ and (\ref{NPexpr5}) to eliminate $\Theta_1$ and $\Theta_2$ the governing evolution  equations (\ref{tory1})-(\ref{tory2}) can be rewritten  in the intuitively appealing form
\be
\frac{\dd}{\dd t}\left(A(x) \pb \right) +  \frac{\dd}{\dd x}\left( A(x) \jpb \right)&=&0, \quad \mbox{where} \quad \jpb= -\kappa_p \pb  \frac{\dd}{\dd x} \left(\mu_e  + \phit \right),\label{Pbcons2} \\ 
\frac{\dd}{\dd t}\left(A(x) \nb \right) +  \frac{\dd}{\dd x}\left( A(x) \jnb \right) &=&0, \quad \mbox{where} \quad  \jnb= -\kappa_n  \nb \frac{\dd}{\dd x} \left(\mu_e  -\phit \right). \label{Nbcons2}
\ee
Furthermore,  $\lambda$ can be expressed in terms of $\mu_e$ as follows
\be
\lambda= \frac{\Lambda}{R} e^{-\mu_e/4}.
\ee
Thus the reformulation of the Quasi-1D PNP model consists of a straightforward method for evaluating the two functional dependence of $G_1(\mu_e; x)$ and $G_2(\mu_e; x)$ on $\mu_e$ and $x$ (contained in (\ref{sinh})-(\ref{lambet}) and (\ref{G1})-(\ref{G2})) and the two coupled parabolic PDEs for $\mu_e$ and $\phit$, (\ref{NPexpr4})-(\ref{Nbcons2}).  A further simplification can be obtained from (\ref{Gdiff}), the relation between $G_1$ and $G_2$, from which we can deduce the local charge neutrality condition
\be
A(x) (\pb -\nb) +\Sigma_l(x)=0, \qquad \Sigma_l(x)= \Lambda^2 \Upsilon  (2 \pi R(x) \sigma(x)). \label{neut}
\ee
Here $\Sigma_l(x)$ represents the fixed charge per unit length (in appropriate dimensionless form) on the wall of the pore. In effect this relation means that we only need to calculate one of the expressions $G_1(\mu_e; x)$ or $G_2(\mu_e; x)$, use this to determine either $\pb$ or $\nb$ from (\ref{NPexpr4}), and evaluate the other from the relation (\ref{neut}).

\paragraph{Calculating the steady state solution} In practice we are usually only interested in the steady state solution to (\ref{NPexpr4})-(\ref{Nbcons2}). Neglecting the time derivatives in (\ref{Pbcons2})-(\ref{Nbcons2}), summing the two equations and taking their difference  yields to the following two equations
\bes
\frac{\dd}{\dd x}\left( A(x) \left( (\nb+\pb)  \frac{\dd\mu_e }{\dd x}- (  \nb-\pb) \frac{\dd \phit }{\dd x}\right) \right)&=&0,\\
\frac{\dd}{\dd x}\left( A(x) \left( (\nb-\pb)  \frac{\dd\mu_e }{\dd x} -(  \nb+\pb) \frac{\dd \phit }{\dd x} \right)\right)&=&0.
\ees
We now write 
\be
\nb+\pb=e^{\mu_e} \Psh(\mu_e;x), \quad \mbox{where} \quad \Psh(\mu_e;x)=G_1+G_2=\int_0^1 \xi \left( e^{\psi} -e^{-\psi} \right) d \xi. \label{Psidef}
\ee
and substitute for $(\nb-\pb)$ from (\ref{neut}) in order to obtain two coupled ODEs for $\mu_e$ and $\phit$
\be
\frac{\dd}{\dd x}\left( A(x) e^{\mu_e} \Psh(\mu_e;x)   \frac{\dd\mu_e }{\dd x}-  \Sigma_l(x)\frac{\dd \phit }{\dd x} \right)&=&0, \label{steady1} \\
\frac{\dd}{\dd x}\left( \Sigma_l(x)  \frac{\dd\mu_e }{\dd x} - A(x) e^{\mu_e} \Psh(\mu_e;x)   \frac{\dd \phit }{\dd x} \right)&=&0.\label{steady2}
\ee
Note that the function $\Psh(\mu_e;x)=G_1+G_2$ can be obtained either by direct solution for $\psi(\xi;\beta,\lambda)$ from (\ref{sinh})-(\ref{sinh_bc}) in which $\lambda= {\Lambda}e^{-\mu_e/4}/R$, or (in the large $\beta$ limit) from the uniformly 
valid asymptotic expression for $G_1$ discussed in Appendix \ref{B} equation (\ref{gsm}) and using the relation (\ref{Gdiff}) to evaluate $G_2$.

%\paragraph{Remarks.} The point about the formulation (\ref{steady1})-(\ref{steady2}) is that it should be the easiest way to deal with steady state problem especially asymptotically. Note in particular that the percentage error that you obtain for the asymptotic expressions for $\Psh$ should be considerably better than for that for $G_2$. It would be good to plot the asymptotic expressions for $\Psh$ versus the fully numerical expressions.
% 
\section{Numerical methods and Results}\label{numerics}
In this section we will present numerical methods for both the full 2D PNP system, the 1D Area Averaged PNP system as well as for the asymptotic Quasi-1D PNP model developed in \S\ref{medium_regime}. 
They will be used to compare the results for two examples: (I) a trumpet shaped pore (see figures \ref{heatmap_2D} \& \ref{heatmap_asym}) and (II) a conical pore geometry (see figures \ref{heatmap_bigpore_2D} \& \ref{heatmap_bigpore_2D_asym}). 

\paragraph{The Quasi-1D PNP solver.}
 The numerical solver of the Quasi-1D PNP  is based on the uniformly-valid large $\beta$ expression for $G_1$ (\ref{gsm}) and on the identity (\ref{Gdiff}), relating $G_2$ to $G_1$. This thus obviates the need to solve the Poisson equation (\ref{sinh})-(\ref{sinh_bc}) for $\psi(\xi,\lambda(x,),\beta(x))$ at every value of $x$.  Instead it only requires 
the solution of the 1D (stationary) continuity equations \eqref{tory1}-\eqref{Theta} in $x$. This represents a very considerable reduction in computational complexity and gives a very fast method, which is particularly suited for the calculation of IV curves. 
Finally, for an applied voltage above a certain threshold, we introduce a relaxation in the iteration. Once $Q$ and $S$ are known, we use \eqref{leadflux} to 
calculate the total current $I$. The full iterative procedure is detailed in Algorithm \ref{alg1}.
\begin{algorithm}[h!]
\SetAlgoLined
 Set $S(x)= p_0(x)/\exp(-\phi_0(x))$,  for $x  \in \{0,\bar l \}$ \; 
 Set $Q(x)= n_0(x)/\exp(\phi_0(x))$,  for $x  \in \{0,\bar l \}$ \;
 Initialise $Q^0(x),S^0(x)$\;
 \While{ $\text{err} > \varepsilon$  \text{ \bf and } $\text{max}\_ {\text{iter}}>m$ }{
 Calculate $\lambda(x)$ using \eqref{lambet}\;
 Calculate $G_1^{m+1} (x)$  and $G_2^{m+1}(x)$ using interpolation between \eqref{G1G2bigbeta}  and \eqref{G1G2smalllambda} \;
  Using $G_1^{m+1} (x)$ and $G_2^{m+1} (x)$  and  equations \eqref{tory1}-\eqref{tory2} calculate $Q^{m+1/2}$, $S^{m+1/2}$ \;
  \uIf{$|V_{appl}|\ge V_c$}{ $Q^{m+1} = \theta Q^{m+1/2} + (1-\theta) Q^{m}$, $S^{m+1} = \theta S^{m+1/2} + (1-\theta) S^{m}$  \;}
    \Else{ $Q^{m+1} = Q^{m+1/2}$, $S^{m+1} = S^{m+1/2}$}
 $\textit{err} =  \|  Q^{m+1} - Q^m\|_2 + \|  S^{m+1} - S^m\|_2$ \;
 $ m= m+1$ \;}
 Calculate $I$ using \eqref{currSQ}
 \caption{Fixed point scheme to calculate $Q$ and $S$ in the steady state.} \label{alg1}
\end{algorithm}

\paragraph{The 2D PNP solver.}
The full 2D steady state PNP system, \ie equations \eqref{nd1}--\eqref{nd8} is solved using a standard $P1$ finite element discretisation and a 
Scharfetter Gummel iteration, \cite{gummel1964self}. For both geometries, we use a non-uniform mesh strongly refined at the charged pore walls in order to 
properly resolve the Debye layers. The meshes are created using Netgen \cite{schoberl1997netgen}, while we use MATLAB to assemble and solve the corresponding discrete systems. 
We use a similar method to solve the 1D Area Averaged PNP system  \eqref{1D_PNPa}- \eqref{1D_PNPc}.

\subsection{Trumpet shaped pores}
\begin{figure}[H]
    \centering
    \begin{subfigure}[b]{0.32\textwidth}
        \includegraphics[width=\textwidth]{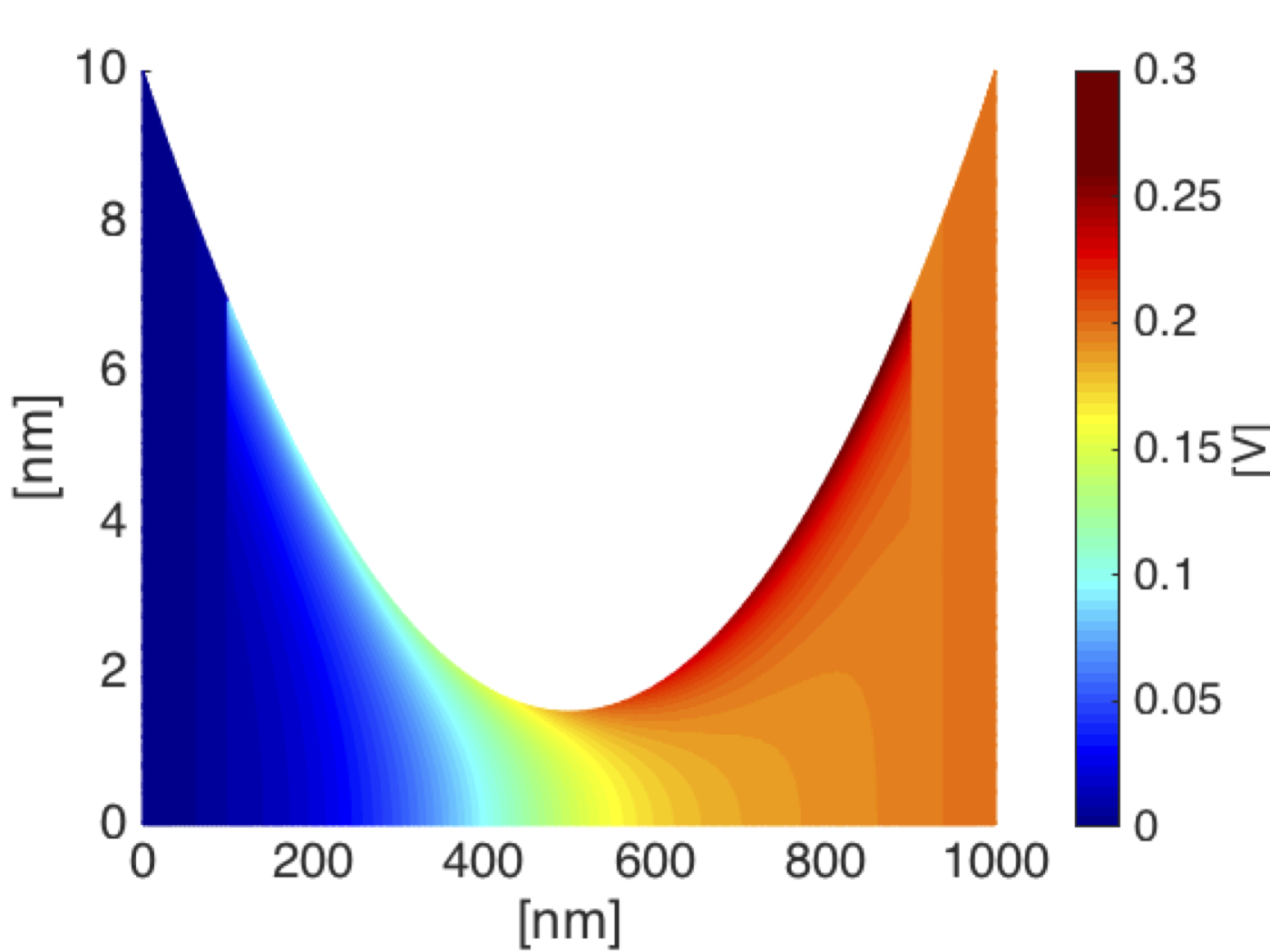}
         \caption{ Potential.  }
    \end{subfigure}
    \begin{subfigure}[b]{0.32\textwidth}
        \includegraphics[width=\textwidth]{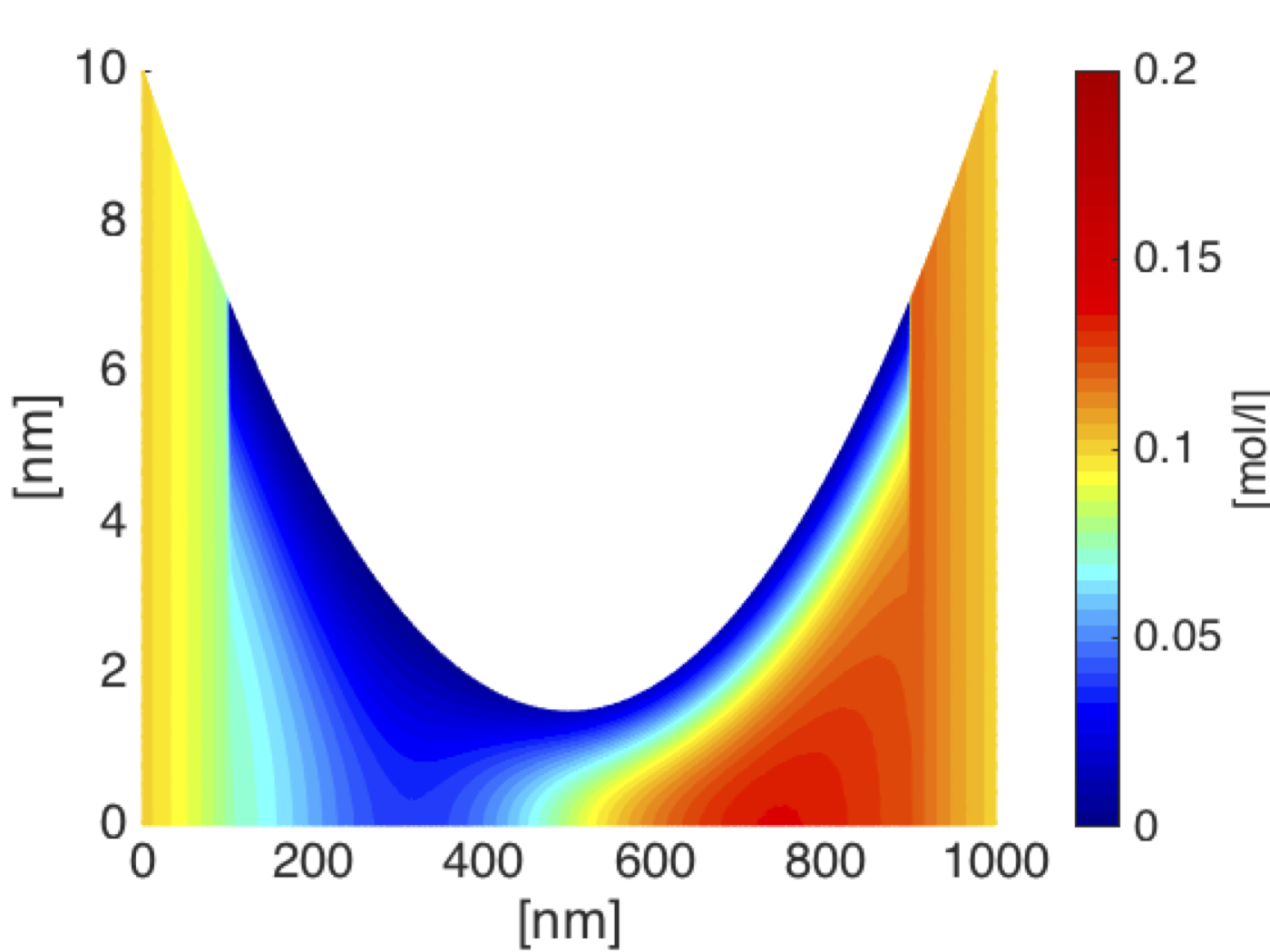}
         \caption{Positive ions conc.}
      \end{subfigure}
      \begin{subfigure}[b]{0.32\textwidth}
         \includegraphics[width=\textwidth]{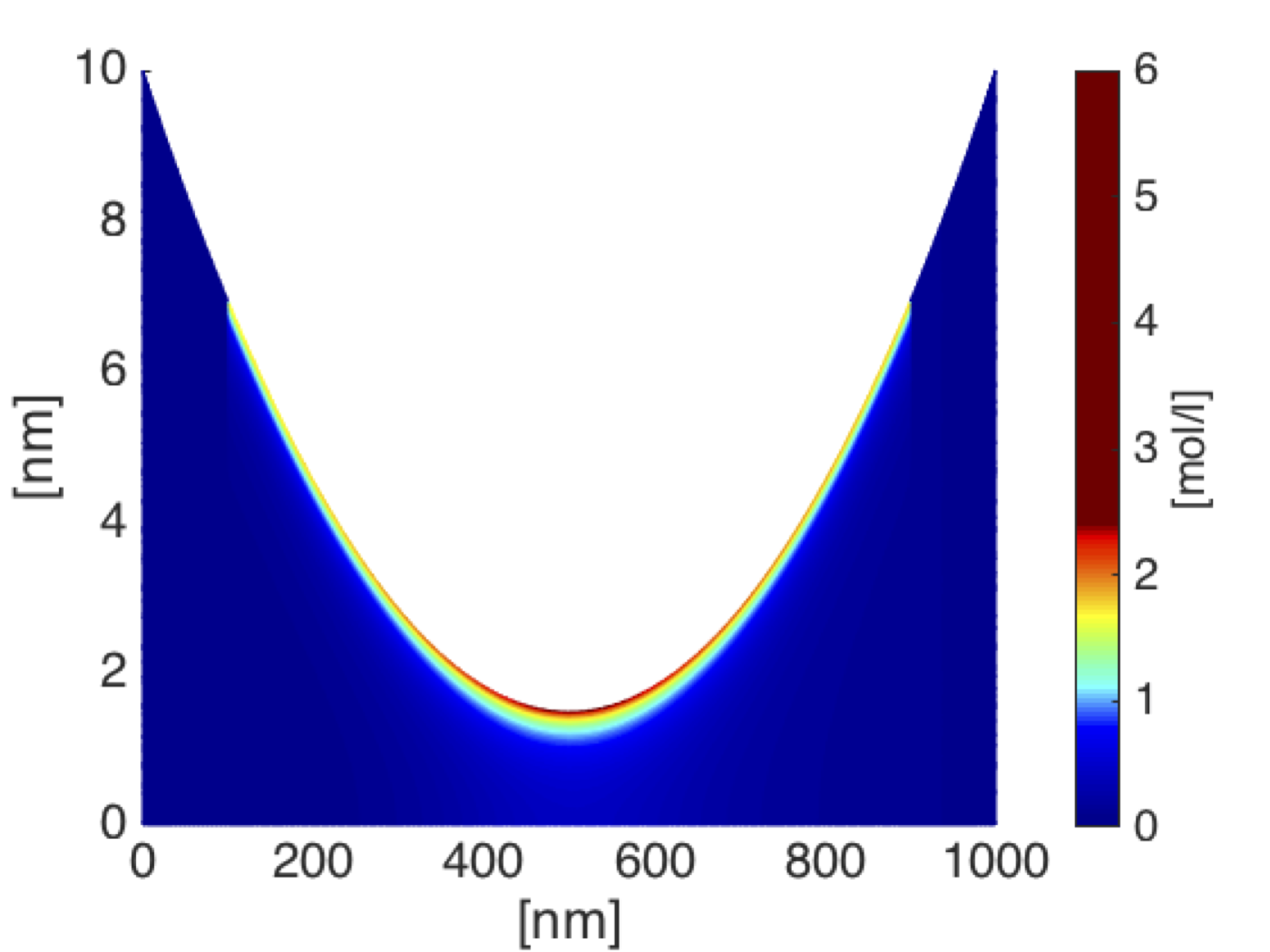}
          \caption{Negative ions  conc.}
      \end{subfigure}
       \caption{ Heat maps of the potential and two ionic concentrations obtained using the 2D PNP solver.     }
       \label{heatmap_2D}
\end{figure}

\begin{figure}[H]
    \centering
    \begin{subfigure}[b]{0.32\textwidth}
        \includegraphics[width=\textwidth]{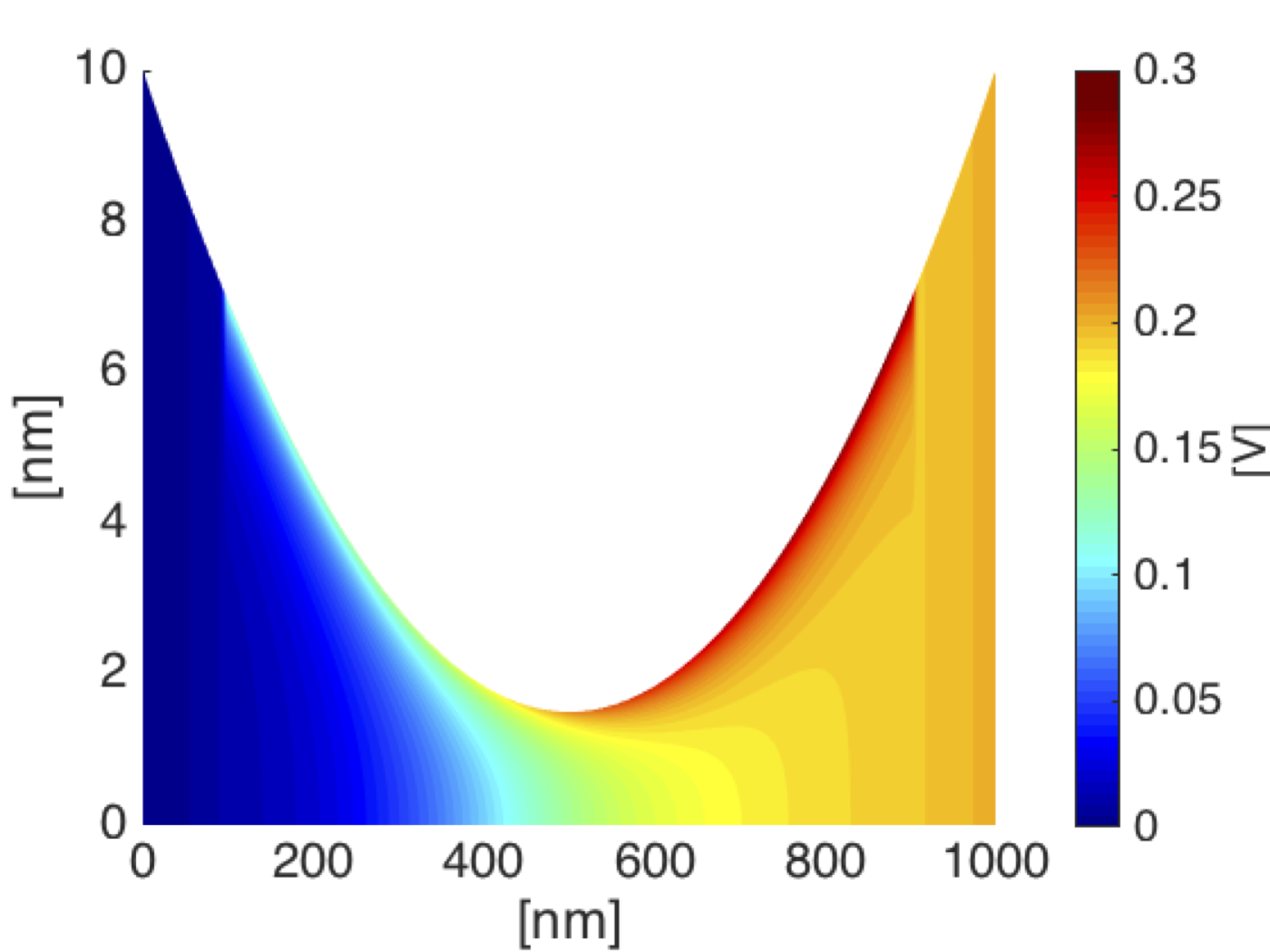}
         \caption{ Potential.  }
    \end{subfigure}
      \begin{subfigure}[b]{0.32\textwidth}
           \includegraphics[width=\textwidth]{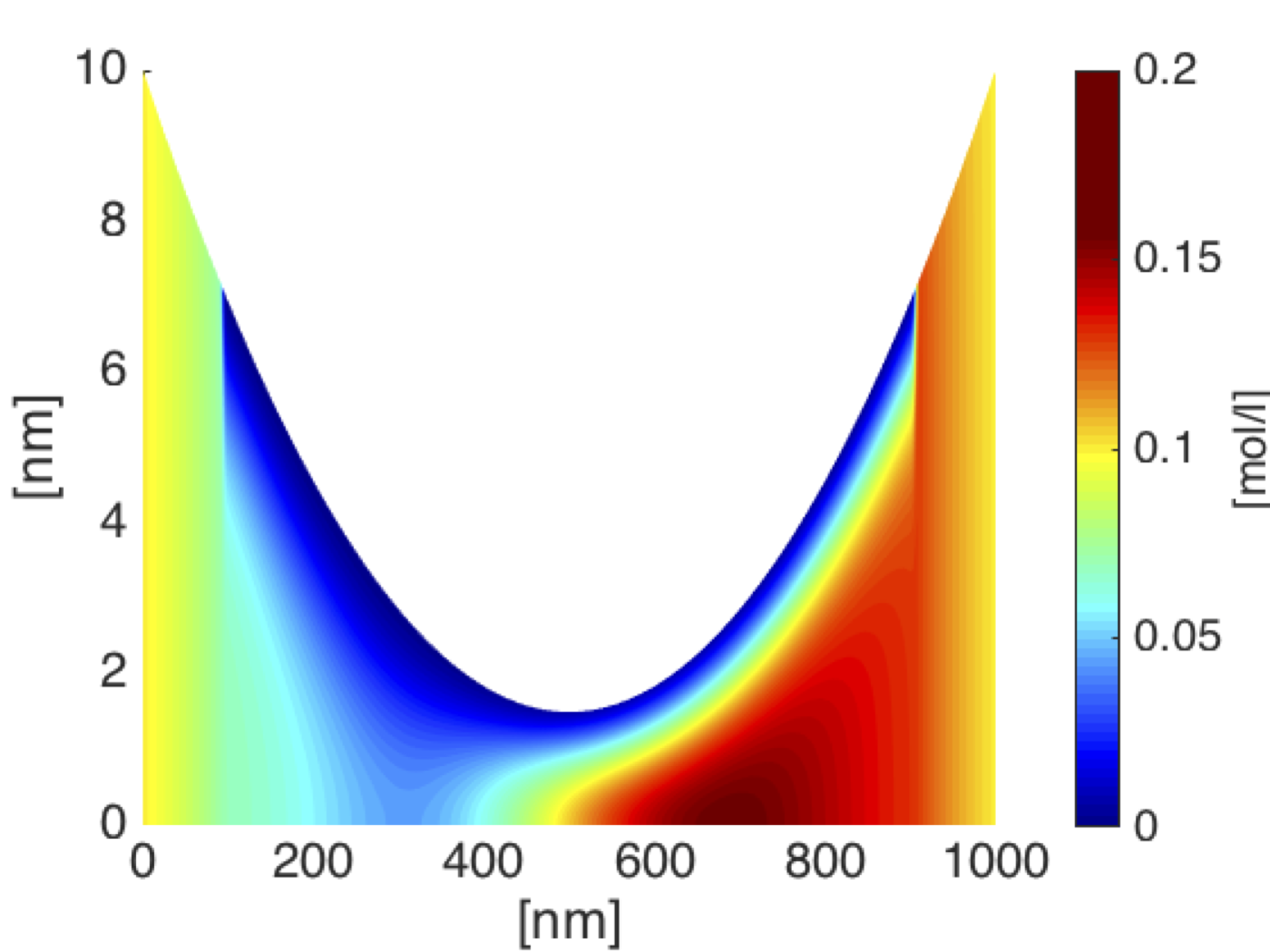}
            \caption{Positive ions conc.}
        \end{subfigure}
        \begin{subfigure}[b]{0.32\textwidth}
            \includegraphics[width=\textwidth]{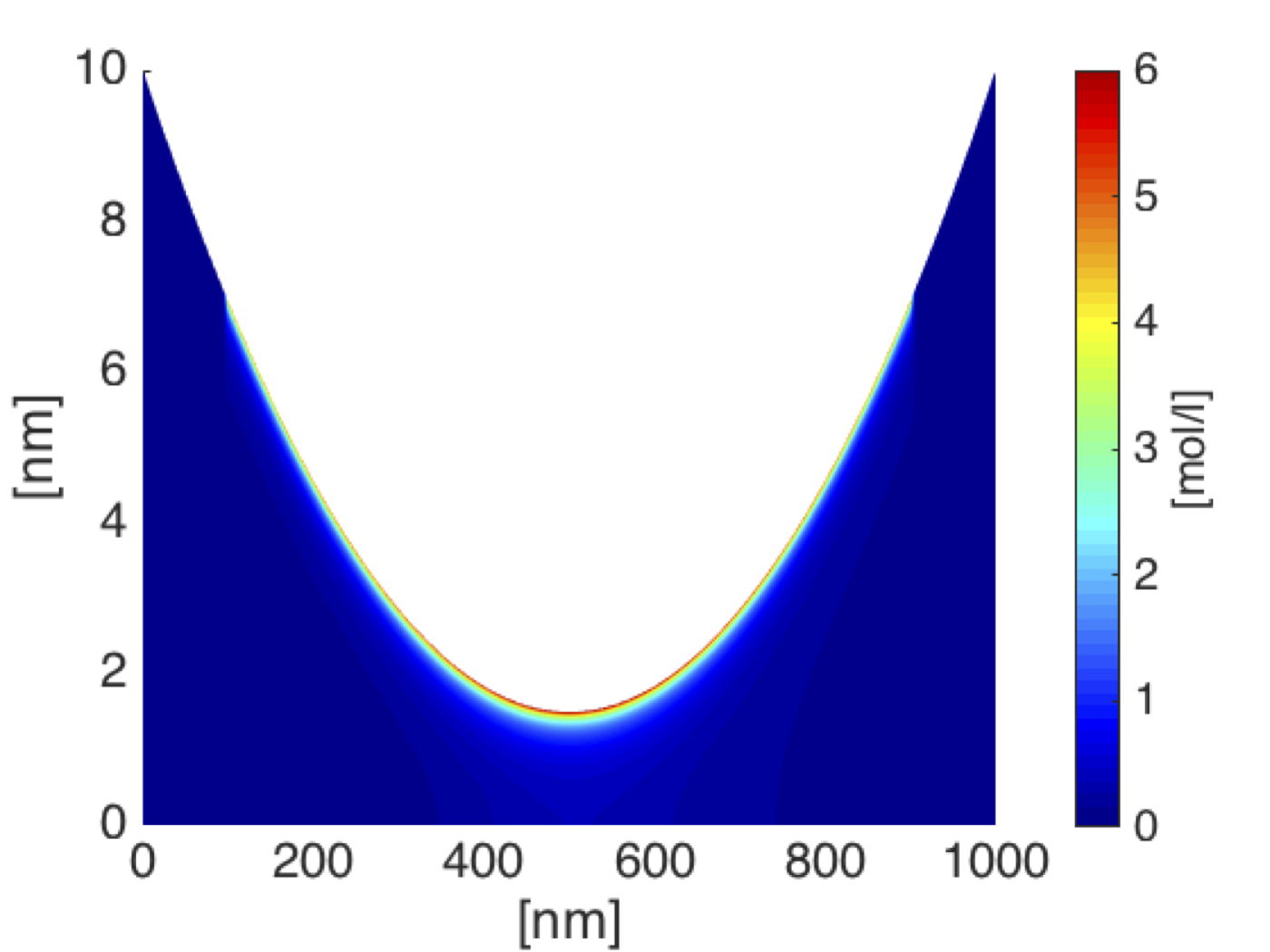}
             \caption{Negative ions  conc.}
      \end{subfigure}
       \caption{Heat maps of the potential and two ionic concentrations obtained using the Quasi-1D PNP solver.}
       \label{heatmap_asym}
\end{figure}

We consider a trumpet shaped pore of length $1000$nm and a radius varying form  $1.5$nm to $10$nm. The corresponding radius is given by
 $ r(x) = 10^{-6}(34x^2 - 34x + 10)$, where both $r$ and $x$ are measured in units of nanometers,  hence the values of $\lambda$ and $\beta$ change continuously with respect to $x$. We set the following parameters:
\be
 V_{appl} = 0.2\text{V and } \quad n_r = n_l = p_r = p_l = 0.1 \text{moles/litre} ~~~~~~\\
 \text{Surface charge profile} \qquad \sigma = \left\{ \begin{array}{cccc} 1 \text{e/nm$^2$} & \text{for} & 100 \text{nm}<x<900\text{nm}\\
0 \text{e/nm$^2$} & \text{for} & |x-500|>400\text{nm}\end{array} \right.  .
 \ee
To obtain accurate and precise  results for the 2D solver a mesh of $360 000$ triangular elements was used. The results of the Quasi-1D and Area Averaged PNP were obtained using a discretization
of $1000$ intervals. Figure \ref{heatmap_2D} and \ref{heatmap_asym} show the solutions to the 2D PNP model and the Quasi-1D PNP model, respectively.

\begin{figure}[H]
     \begin{center}
       \includegraphics[width=1\textwidth]{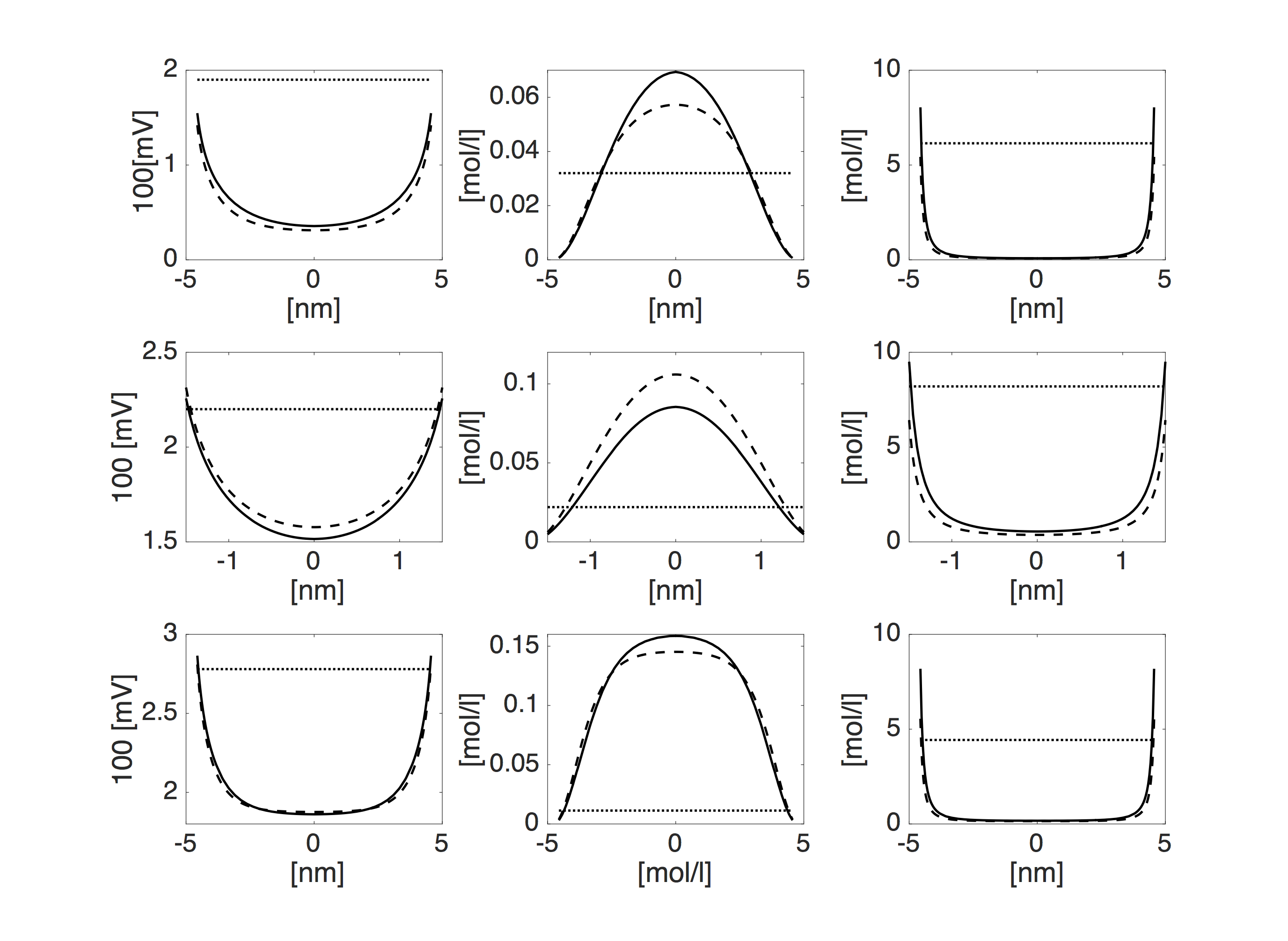}
        \caption{Comparison of the potential (left column), positive ions concentrations (centre column) and negative ions concentrations (right column) calculated 
        over the cross-section at $x = 200$nm (top row), $x = 500$nm (middle row) and $x = 800$nm (bottom row), obtained using the 2D finite element solver (solid lines), the 
        1D Area Averaged PNP (dotted lines) and the Quasi-1D PNP solver described in Algorithm \ref{alg1} (dashed lines) for a trumpet shaped shape pore of length $1000$nm and radius varying from $1.5$ to $10$nm.}  
       \label{profiles}
        \end{center}
\end{figure}

In order to compare the results from the two different methods we plot the cross sectional profiles of the potential and concentrations at $ x = 200, 500$ and $800$nm in 
Figure \ref{profiles}. We observe that the solution to the Quasi-1D PNP model is a very good match to that of the full 2D PNP equations. This is especially so for the potential (left column) and the negative ions concentration (right column) for which both solutions have almost identical behaviour within the Debye layers. 

\begin{figure}[H]
      \begin{subfigure}[b]{0.48\textwidth}
        \includegraphics[width=\textwidth]{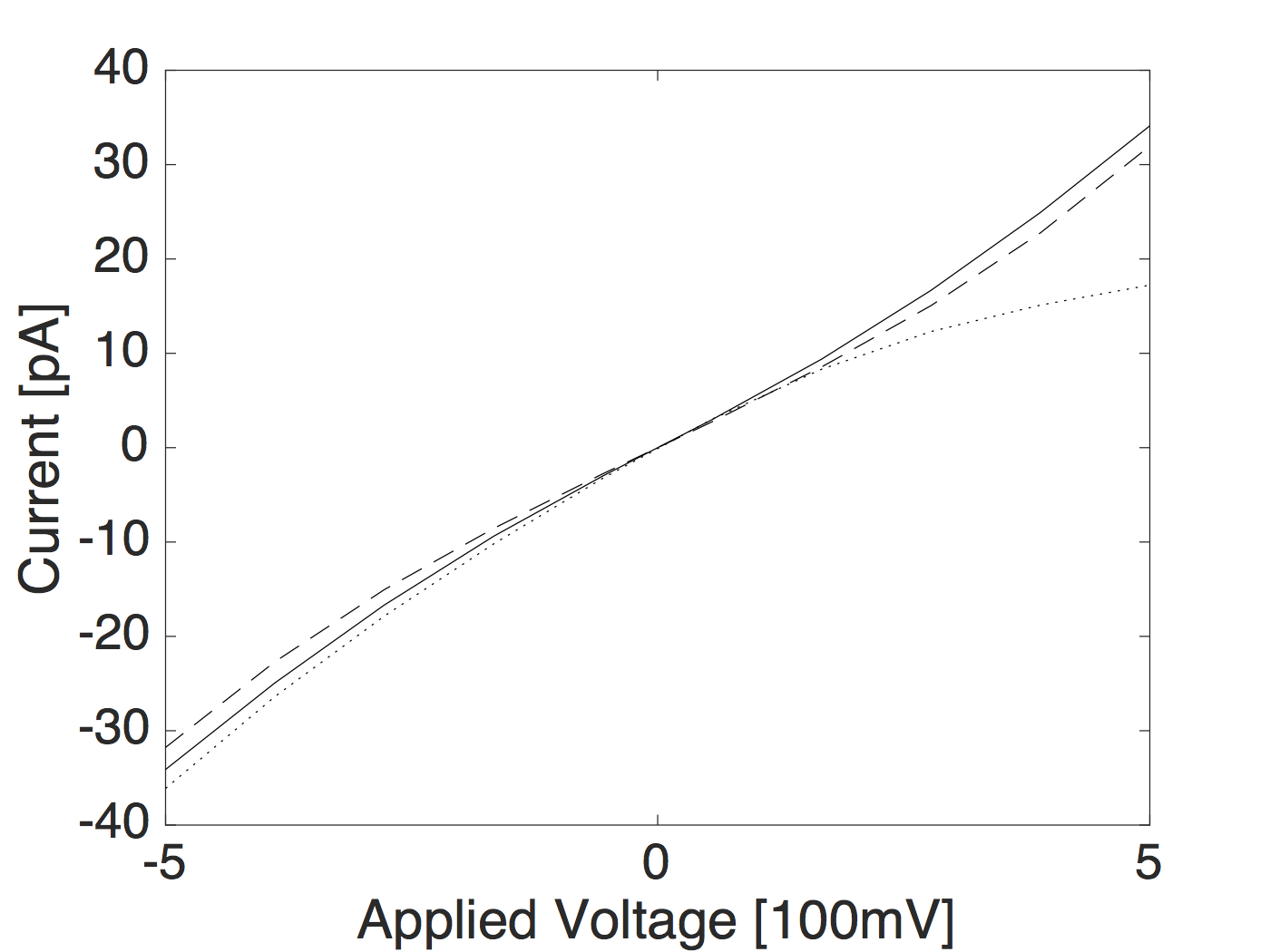}
        \caption{$\sigma$ = 1 e/nm$^2$}
    \end{subfigure}
   \begin{subfigure}[b]{0.48\textwidth}
        \includegraphics[width=\textwidth]{IV_trump_sig02.png}
        \caption{$\sigma$ = 0.2 e/nm$^2$}
    \end{subfigure}	
       \caption{IV curves for surface charges  $\sigma = 0.2$ e/nm$^2$  (right plot) and $\sigma = 1$ e/nm$^2$ (left plot)   obtained using  the Quasi-1D PNP  solver (dashed lines), the 2D PNP solver (straight lines) and 1D Area Averaged PNP method (dotted lines). }
        \label{Iv}
       \end{figure}

Next we compare the IV curves in the case of different surface charge densities $\sigma = 1$ e/nm$^2$ and $\sigma = 0.2$ e/nm$^2$ within the central region of the pore  $|x-500|<400\text{nm}$ (we take $\sigma=0$ outside this region) see Figure \ref{Iv}. We observe very good agreement between results from the full 2D PNP model and the Quasi-1D PNP model for both values of the surface charge density. Notably the agreement of the Area Averaged PNP to the full 2D PNP model is much worse than that of the Quasi-1D PNP model.

% Big pore section 
\subsection{Conical Shaped Pores}
\begin{figure}[H]
    \centering
    \begin{subfigure}[b]{0.32\textwidth}
        \includegraphics[width=\textwidth]{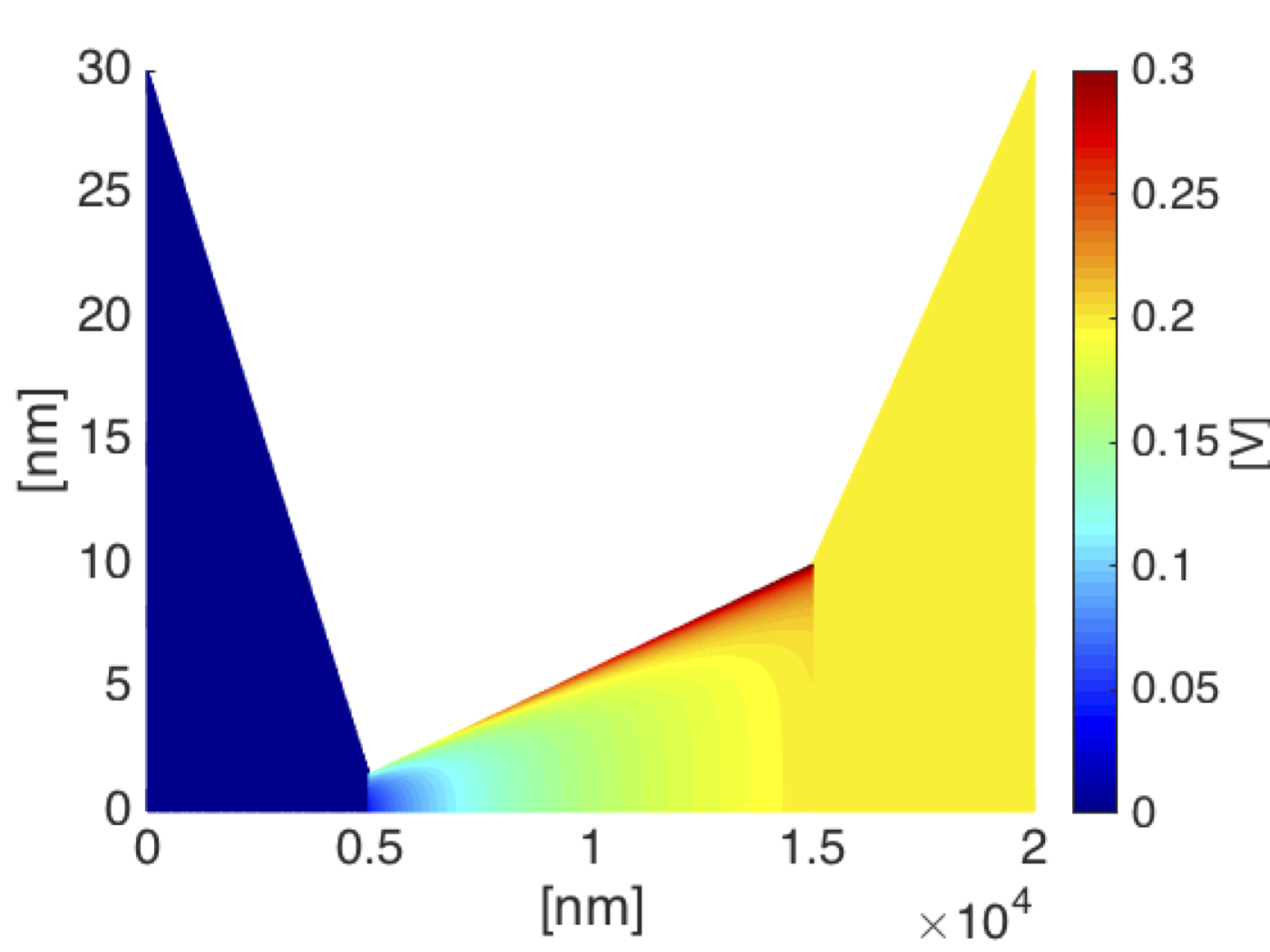}
         \caption{ Potential.  }
    \end{subfigure}
      \begin{subfigure}[b]{0.32\textwidth}
           \includegraphics[width=\textwidth]{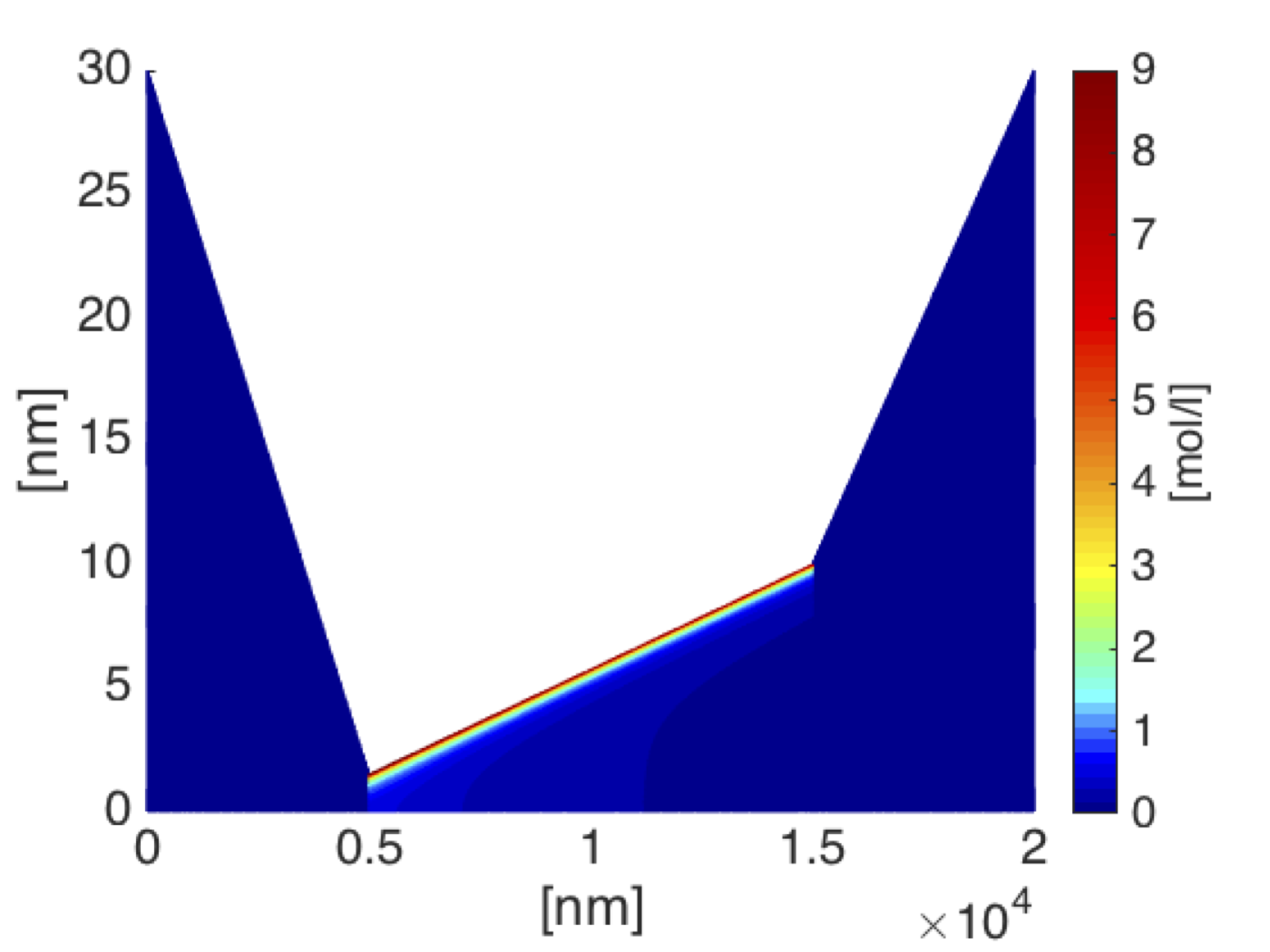}
            \caption{Negative ions conc.}
        \end{subfigure}
        \begin{subfigure}[b]{0.32\textwidth}
            \includegraphics[width=\textwidth]{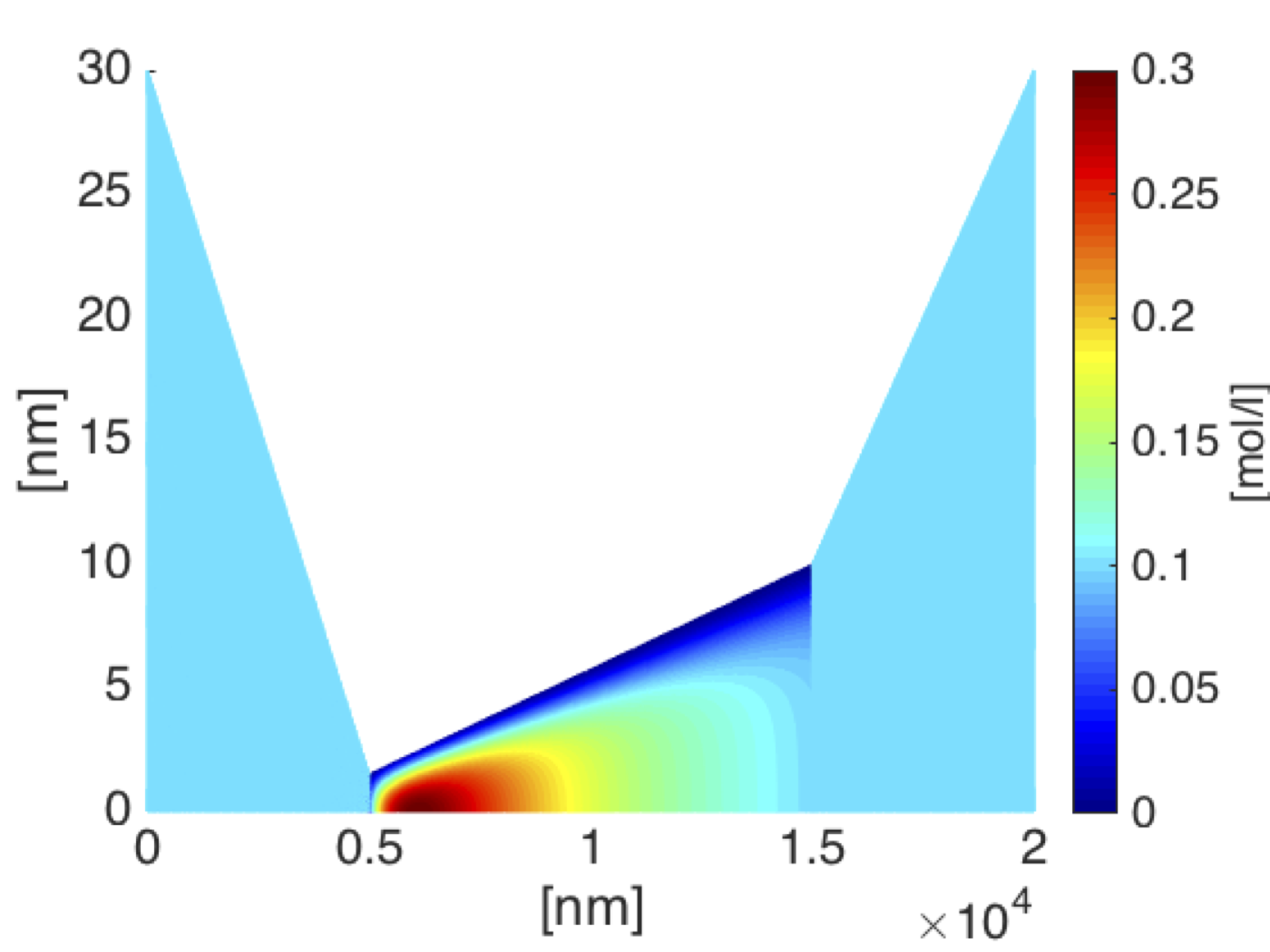}
             \caption{Positive ions  conc.}
      \end{subfigure}
       \caption{ Heat maps of the potential and two ionic concentrations obtained using the 2D PNP solver for the conical pore.  }
       \label{heatmap_bigpore_2D}
\end{figure}

\begin{figure}[H]
    \centering
    \begin{subfigure}[b]{0.32\textwidth}
        \includegraphics[width=\textwidth]{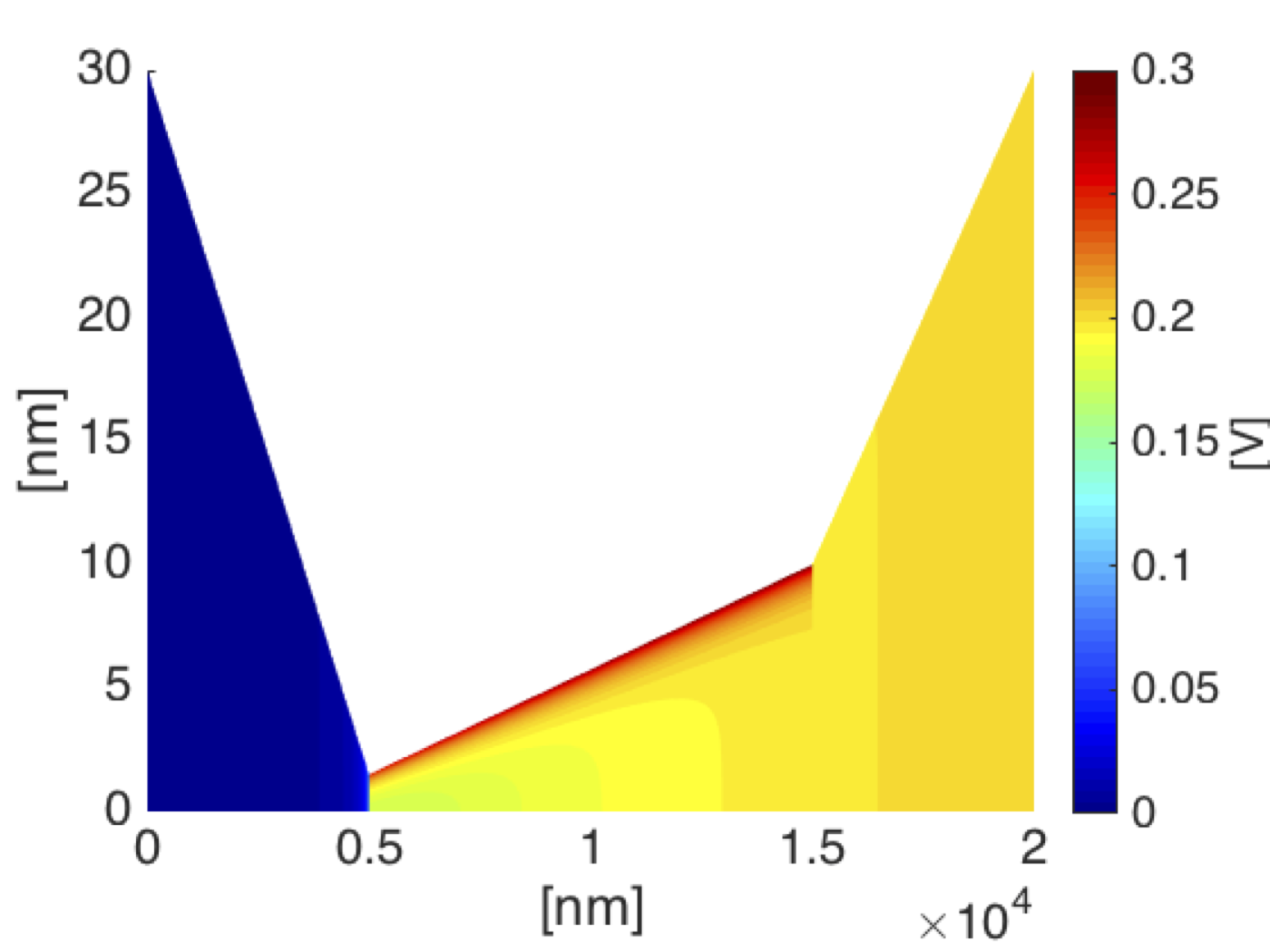}
         \caption{ Potential.  }
    \end{subfigure}
      \begin{subfigure}[b]{0.32\textwidth}
           \includegraphics[width=\textwidth]{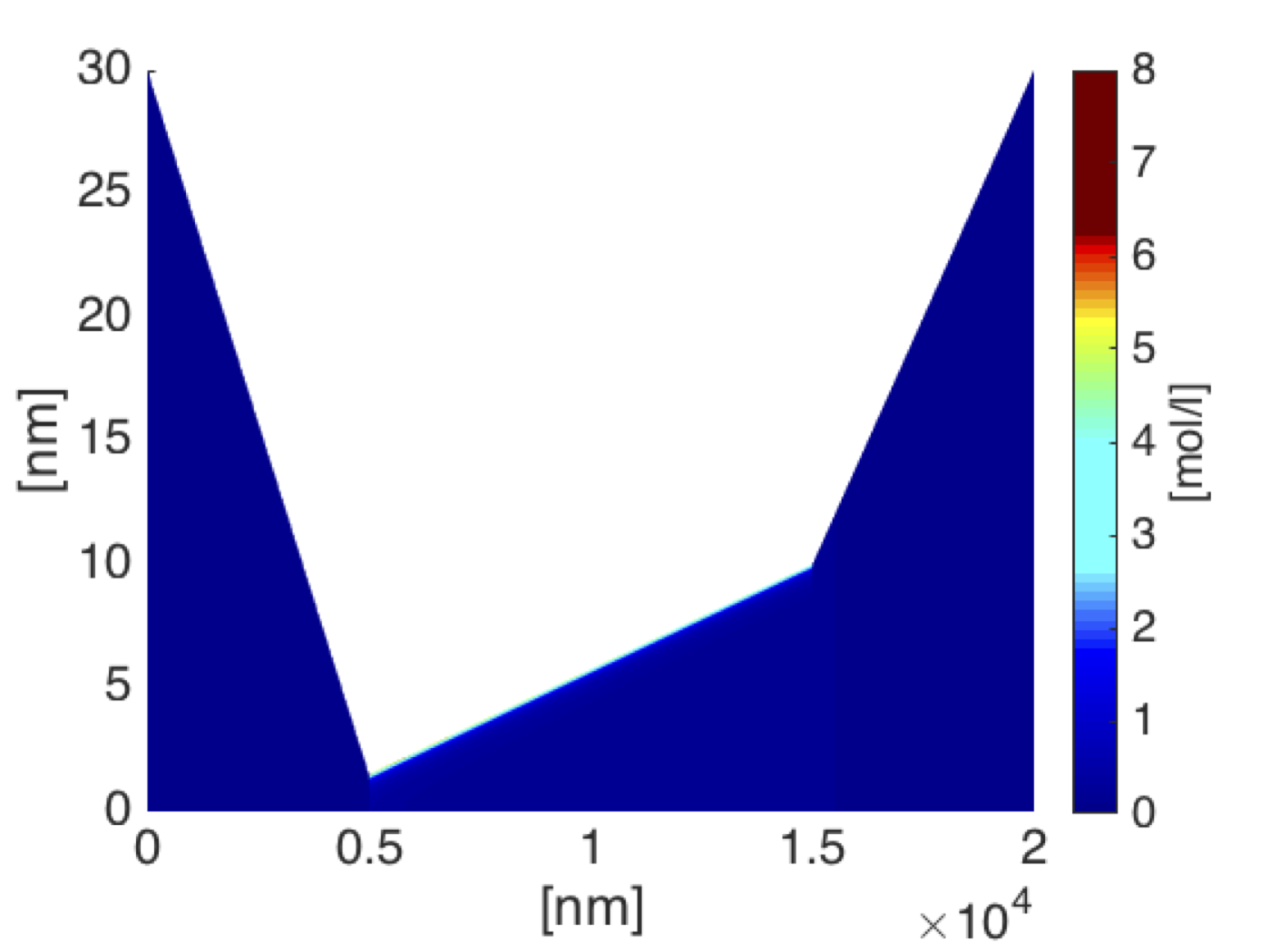}
            \caption{Negative ions conc.}
        \end{subfigure}
        \begin{subfigure}[b]{0.32\textwidth}
            \includegraphics[width=\textwidth]{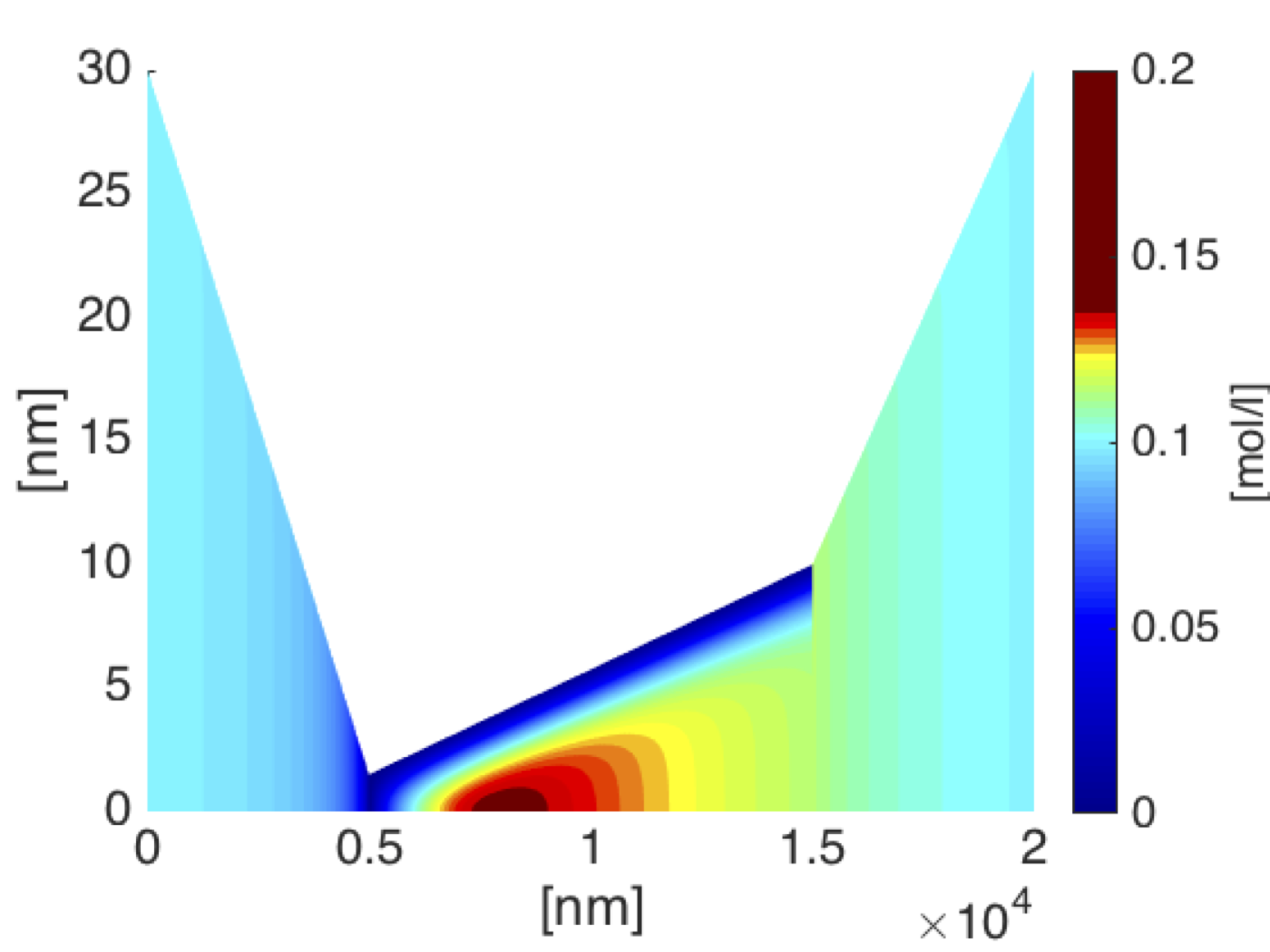}
             \caption{Positive ions  conc.}
      \end{subfigure}
       \caption{ Heat maps of the potential and two ionic concentrations obtained using Quasi-1D PNP solver for the conical pore. }
       \label{heatmap_bigpore_2D_asym}
\end{figure}

\begin{figure}[H]
     \begin{center}
       \includegraphics[width=1\textwidth]{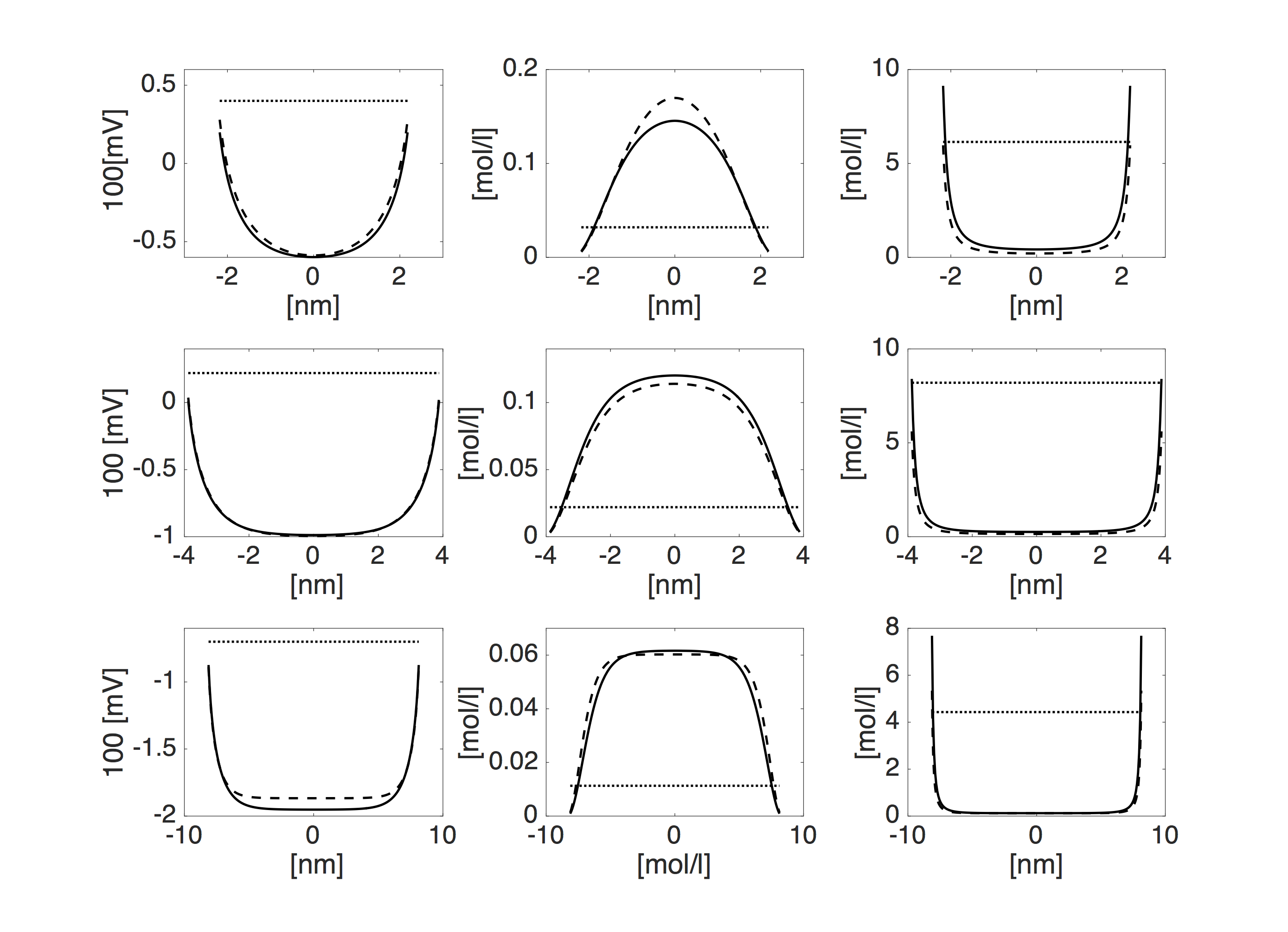}
        \caption{The conical pore. Comparison of the potential $\phi$ (left column), negative ion concentration $n$ (centre column) and positive ion concentration $p$ (right column) calculated over 
        the cross-section at $x = 5800$nm (top row), $x = 7800$nm (middle row) and $x = 12800$nm (bottom row), obtained using the 2D finite element solver (solid lines)
        the 1D Area Averaged PNP (dotted lines) and the Quasi-1D PNP solver described in Algorithm \ref{alg1} (dashed lines) for a linear pore of length $10000$nm and radius varying from $1.5$ to $10$nm.     }  
       \label{profiles_bigpore}
        \end{center}
\end{figure}

Here, motivated by experimental work on etched pores with conical shape \cite{siwy_ss_2003}, we consider 
 a conically shaped pore of length $ 10 000 $nm with radius varying between $1.5$nm and $10$nm  (see figures \ref{heatmap_bigpore_2D} \& \ref{heatmap_bigpore_2D_asym}). This very narrow pore tip is a good model for the tip of a polyethylene terephthalate (PET) nanopore, as used by Siwy\etal \cite{siwy_ss_2003}. 
 It is well known that such narrow 
 tips strongly influences the ion transport through the pore \cite{pietschmann_pccp_2013}. We include two bath regions of 5 $\mu$m length each. Here we consider a pore with uniform surface charge density  inside the pore which corresponds with  5000nm$<x<$15000nm and zero outside this section.
Because of the different length scales and the boundary layer scale we use a highly anisotropic mesh of $ 7 \times 10^5$  triangular elements 
(calculated using Netgen \cite{schoberl1997netgen}) refined at the boundary to capture the boundary effects. 
Figures \ref{heatmap_bigpore_2D} and \ref{heatmap_bigpore_2D_asym} show the results of the full 2D model and those of the Quasi-1D PNP model, respectively. The corresponding cross sectional profiles are depicted in Figure \ref{profiles_bigpore}.

Again we observe very good agreement between the Quasi-1D PNP model solution and the full 2D results close to the charged pore walls. While the  discrepencies between the potentials and the negative ions calculated using these two methods are negligible those for the positive ion concentrations are more marked. 

\begin{figure}[H]
      \begin{subfigure}[b]{0.48\textwidth}
        \includegraphics[width=\textwidth]{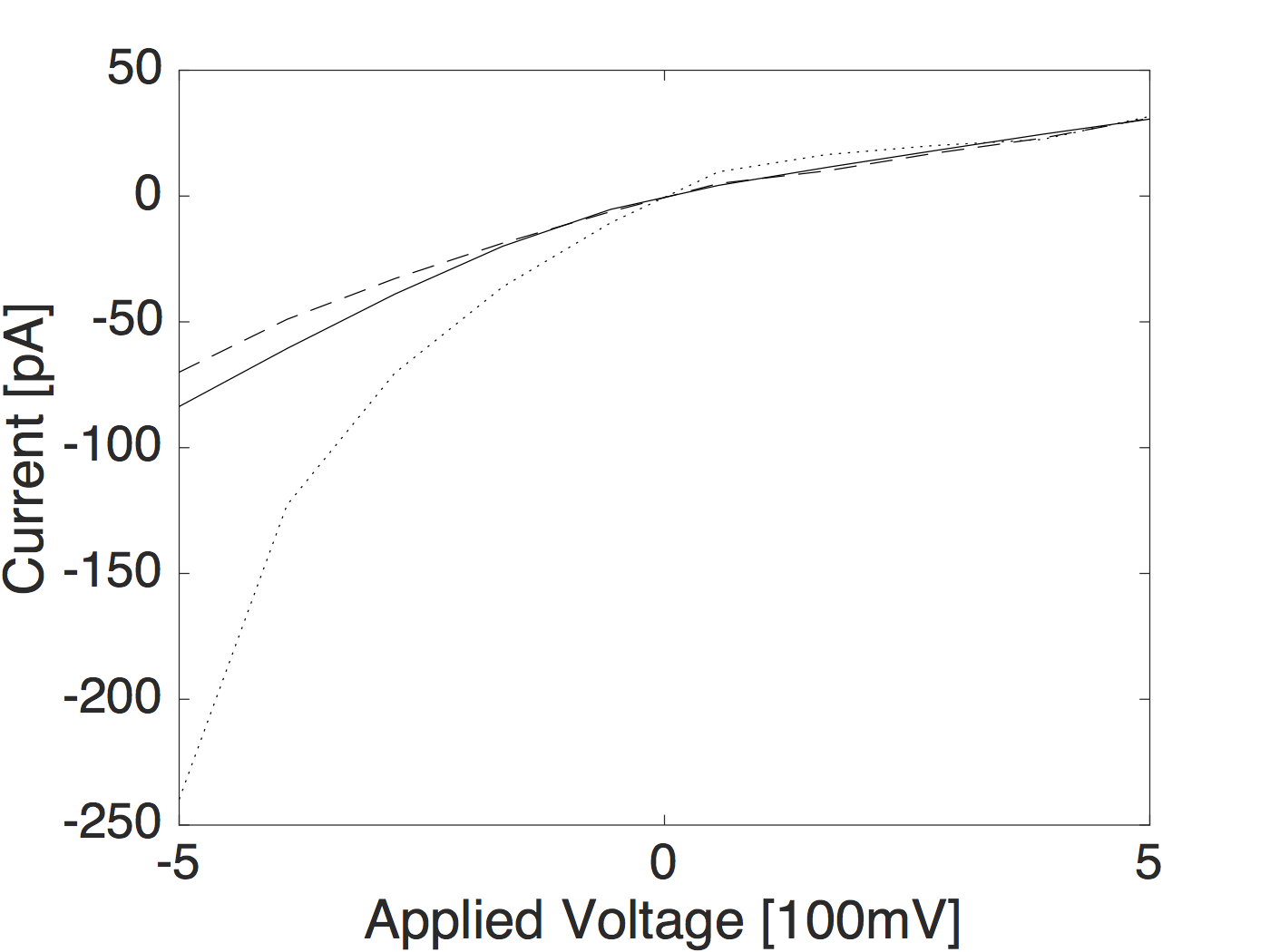}
        \caption{$\sigma$ = 1 e/nm$^2$}
    \end{subfigure}
   \begin{subfigure}[b]{0.48\textwidth}
        \includegraphics[width=\textwidth]{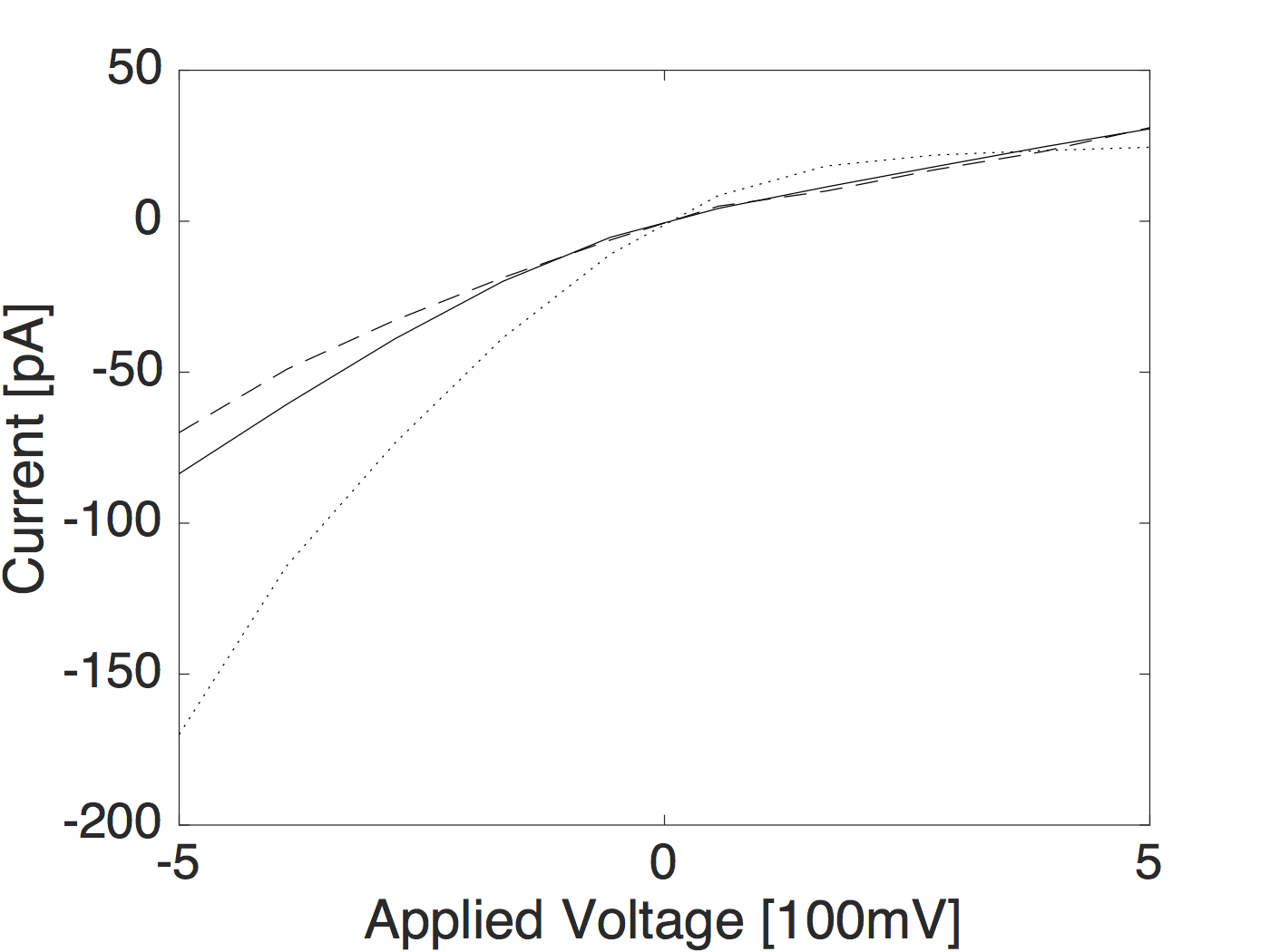}
        \caption{$\sigma$ = 0.2 e/nm$^2$}
    \end{subfigure}	
   \caption{IV curves  for the conical pore obtained using the surface charges  $\sigma = 0.2$ e/nm$^2$   (right plot)  and $\sigma = 1$ e/nm$^2$   (left plot), in the region 5000nm $<x<$ 15000nm, obtained using  the  Quasi-1D PNP (dashed lines), 2D PNP solver (solid lines) and 1D Area Averaged PNP method (dotted lines).}
   \label{bigpore_IV}
 \end{figure}

Finally Figure \ref{bigpore_IV} shows the IV curves obtained from the 2D FEM code, the Quasi-1D PNP solver and the 
Area Averaged PNP equations.  
There is much better agreement between the full 2D solver and the Quasi-1D PNP solver than between either of these and the 1D Area Averaged PNP solver (this again overestimates the influence of the  geometrical asymmetry of the pore and surface charge influence  on the current). Note that the 
Quasi-1D PNP solver captures the nonlinear IV curve, and the corresponding rectification behaviour, much better than the 1D Area Averaged PNP solver.

 \section{Conclusion}
 
 In this work we applied asymptotic methods to a two dimensional Poisson-Nernst-Planck  (PNP) model, for the transport processes occurring within a long thin electrolyte filled nanopore with charged walls, in order to systematically derive a reduced order model for ion transport within the nanopore. We term this the Quasi-1D PNP model. In  order to investigate the validity of this novel model we conducted numerical experiments on two different nanopore geometries in which we compared results from the Quasi-1D PNP model to solutions of the full two dimensional PNP model, which we solved using a finite element method.  In the geometries we considered the comparison between the two approaches was very favourable and furthermore the computational cost of solving the reduced order model was many times less than that for solution of the full 2D model, which requires the use of a very large number of finite elements in order to obtain sufficient accuracy. In addition, we also compared the solution of these two models to the solutions of the one-dimensional  Area Averaged PNP  equations, which is a commonly used approximation of the PNP model in nanopores, and showed that this model gives a poor representation of the full PNP  equations. In this context we also note that the Area Averaged PNP equations are also widely applied to biological ion channels \cite{singer_gillespie_norbury_eisenberg_2008,singer2009poisson,ionreport} but that no comparison has yet been made between numerical solutions to the PNP equations in 3D and solutions to the 1D Area Averaged equations in an ion channel geometry. Furthermore, given that the Debye length in intra- and extra-cellular fluid ($\approx$0.14 Molar) is around 1.3nm, and that the narrow neck of an ion channel is around 0.4nm (comparable to the Debye length), one might expect that the Quasi-1D PNP provides at least as good an approximation (if not better) to the full 3D PNP as the Area Averaged PNP (which should only be valid if the dimensions of the channels are much smaller than the Debye Length).

The numerical experiments presented here confirm the validity of the assumptions made in the derivation of the Quasi-1D PNP equations. We observe that the method resolves the 
behaviour of solutions inside the Debye layers correctly and gives substantially better results then the commonly used 1D area averaged approximations. 
Since surface charge influences the transportation and rectification behaviour of the pore significantly, the correct
resolution of the numerical simulations is of great importance. The proposed asymptotics serves as a starting point for further developments in this direction, in particular
\begin{itemize}
 \item the efficient implementation of a 1D solver to calculate IV curves for nanopores
 \item the extension of the asymptotic analysis for nonlinear PNP models
 \item and the comparison on the results with experimental data.
\end{itemize}

\section*{Acknowledgements} The work of JFP was supported by DFG via Grant 1073/1-2. MTW and BM acknowledges financial support from the Austrian Academy of Sciences \"OAW via the New Frontiers 
Grant NST-001. BM acknowledges the support from the National Science Center from award No DEC-2013/09/D/ST1/03692.

\bibliographystyle{plain}
\bibliography{references}
\appendix

\section{Interpolating the function $G_1(\lambda,\beta)$}\label{B}
In trying to obtain approximate expressions for $G_1$ and $G_2$ we made use of the fact that $\beta \gg 1$ for most practical applications of interest. With this proviso we were able to obtain the behaviour of these functions for both $\lambda=O(1)$ and $\lambda \ll 1$, in (\ref{G1G2bigbeta}) and (\ref{G1G2smalllambda}), respectively. In the case of $G_2$ there is enough overlap between the $\lambda=O(1)$ and $\lambda \ll 1$ limits to obtain a good approximation of this function for all values of $\lambda$  (provided $\beta \gg 1$). Furthermore if we can obtain a good approximation to $G_1$ we can use the exact expression (\ref{Gdiff}) to evaluate $G_2$.

\begin{figure}[H]
      \begin{subfigure}[b]{0.48\textwidth}
        \includegraphics[width=\textwidth]{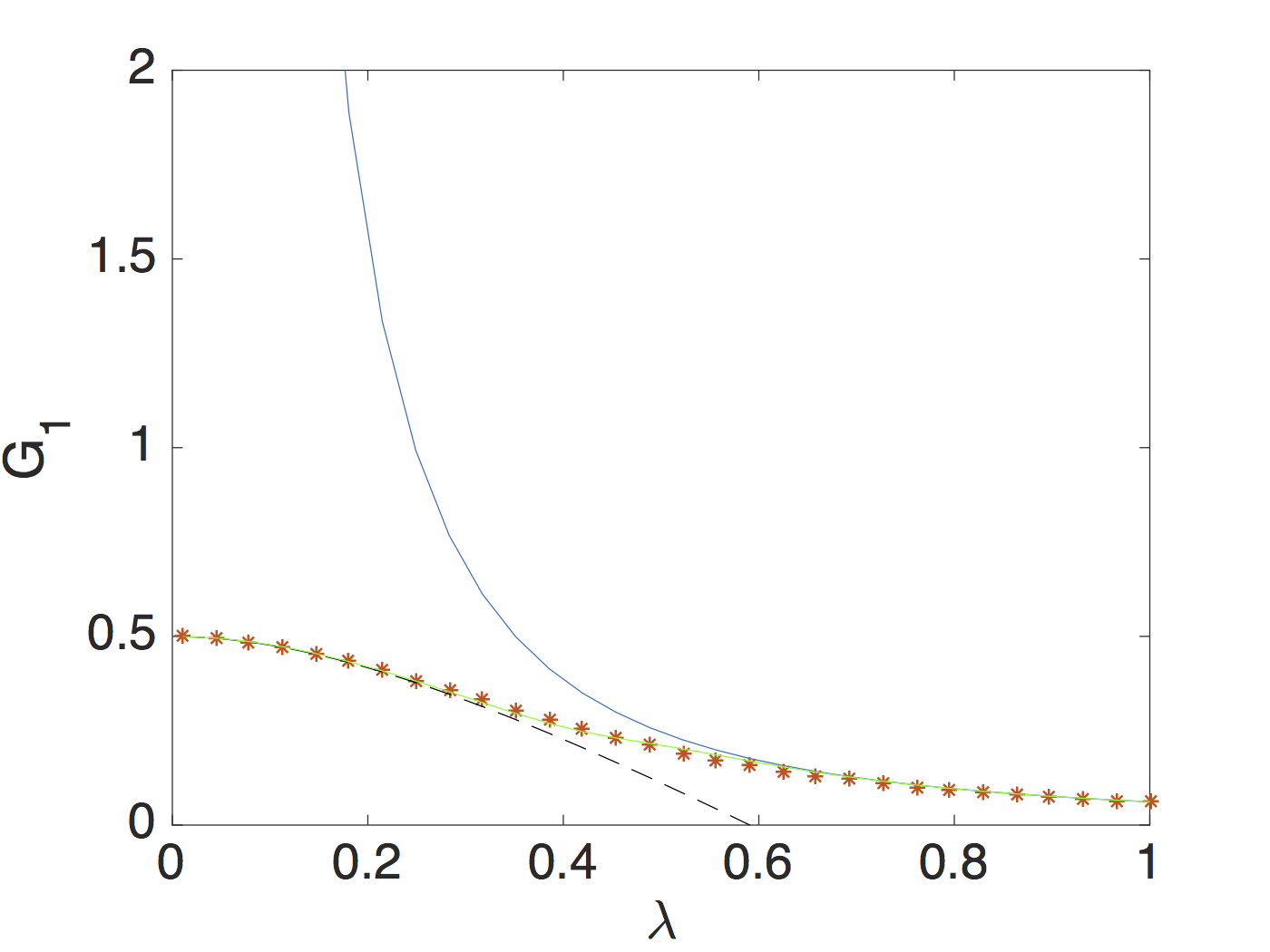}
                \caption{$\beta=5$}
    \end{subfigure}
   \begin{subfigure}[b]{0.48\textwidth}
        \includegraphics[width=\textwidth]{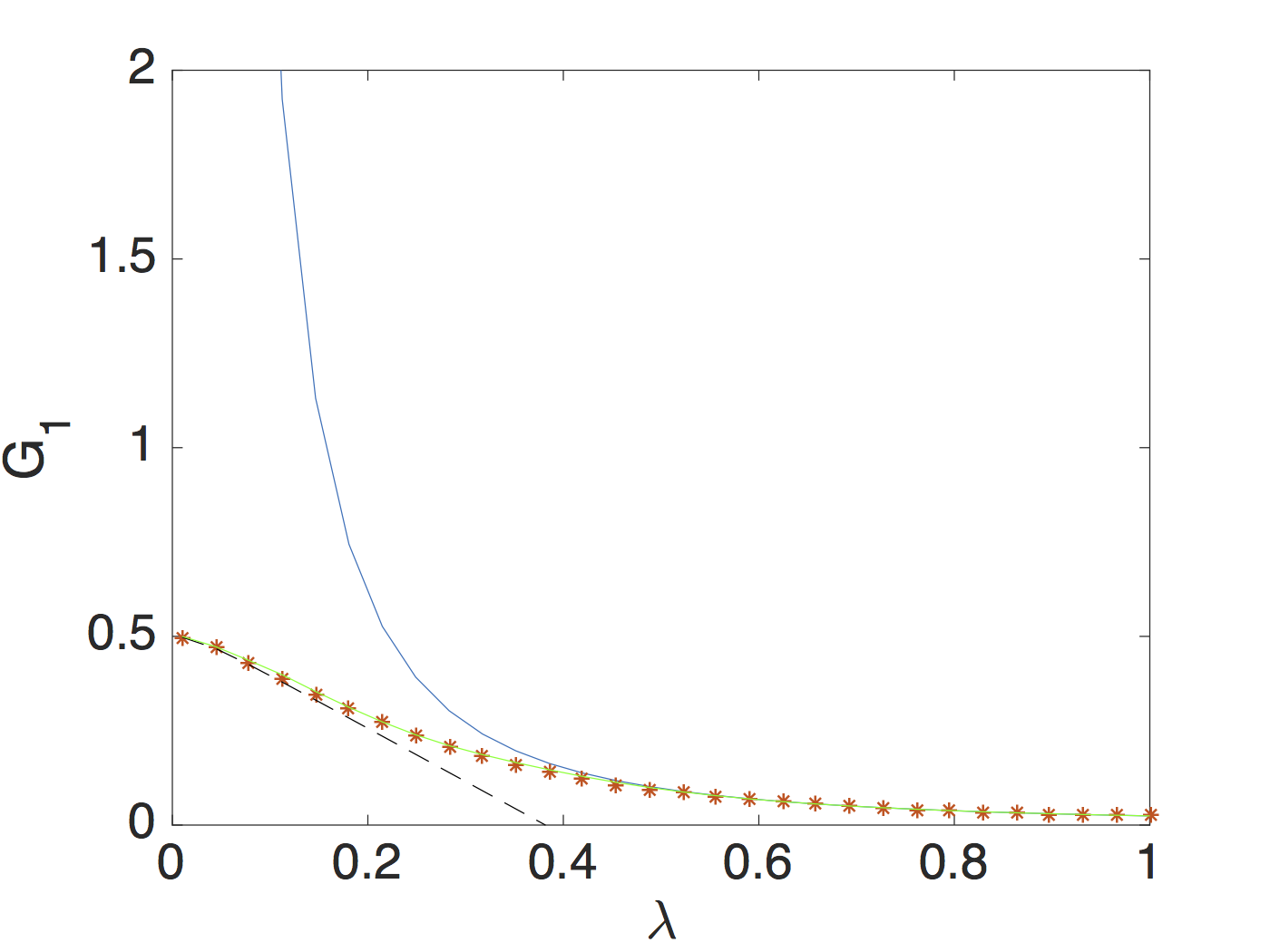}
        \caption{$\beta=50$}
    \end{subfigure}

        \caption{Evaluation of $G_1$ for  $\lambda \ll$  (black line),  $\lambda=O(1)$,   interpolated with described procedure(green line) and exact solution (red dots) for two different values of the parameter $\beta$.}  
          \label{G1_smooth}
\end{figure}

However this is not true of $G_1$ for which there is rapid switching between the $\lambda=O(1)$ and $\lambda \ll 1$ behaviour. In order to obtain a very good approximation of $G_1$ for all values of $\lambda$ we seek to interpolate between the two behaviours by formulating a uniformly valid asymptotic solution. We start by denoting the $O(1)$ and small $\lambda$ behaviours of $G_1$ (as given in  (\ref{G1G2bigbeta}) and (\ref{G1G2smalllambda})) by $\gil$ and $\gis$, noting that they are given by the following functions of $\lambda$ and $\beta$:
\be
\gil(\lambda,\beta)=\frac{1}{48 \lambda^2} \frac{\beta^2+12 \beta+48}{\beta(\beta+4)}, \qquad \gis(\lambda,\beta)=\frac{1}{2} -\frac{ 2 \sqrt{2} \lambda^2 \beta}{\sqrt{8 + \lambda^2 \beta^2} +2 \sqrt{2} +\lambda \beta}. \label{G1SL}
\ee
We now introduce the switching function $S_w(\lambda,\beta)$, which we design to switch smoothly between the two behaviours around some optimal value $\lambda$ denoted by $\lambda_{sw}(\beta)$, this is defined by
\be
S_w(\lambda,\beta)=\frac{1}{2} \left( 1+ \tanh( 12 (\lambda-\lambda_{sw}(\beta) )) \right).
\ee
By fitting to data we find that the optimal switching value is well-approximated by
\be
\lambda_{sw}(\beta) = 0.276+0.9 \beta^{-1}.
\ee
We take $\gsm$, the smoothed approximation of $G_1$ (the uniformly asymptotic solution), to be given by
\be
\gsm(\lambda,\beta)=S_w(\lambda,\beta) \left[ \gil(\lambda,\beta) H(\lambda-0.1) + \gil(0.1,\beta) (1-H(\lambda-0.1)) \right]~~~~~~~~~ \non \\
+ (1-S_w(\lambda,\beta))\gis(\lambda,\beta).~~~~~~~~~~~~~~ \label{gsm}
\ee
Note that we have cutoff the singular behaviour of $\gil(\lambda,\beta)$ as $\lambda \ra 0$ with the use of the Heaviside function $H(\lambda-0.1)$ within the square brackets. Plots of $\gsm(\lambda,\beta)$ against $\lambda$ are made for various values of $\beta$ and compared to the full numerical solution of $G_1$ in figure \ref{G1_smooth}. It can be seen that this uniformly valid asymptotic approximation to $G_1$ is extremely accurate for large $\beta$.

\end{document}